\documentclass[12pt]{article}

\usepackage{amsfonts,amsmath,amssymb,latexsym,xcolor,mathrsfs,tikz,breqn,extarrows}
\usepackage{amsthm,authblk}
\usepackage{multirow}
\usepackage{amssymb}
\usepackage{amsmath}
\usepackage{amsthm}
\usepackage{tikz}
\usepackage{mathtools}
\usepackage{calc}
\usepackage{tikz}
\usepackage{pgfplots}
\usetikzlibrary{patterns,arrows,decorations.pathreplacing}
\usepackage{diagbox}
\usepackage{makecell}
\usepackage{bbding}

\def\q{\hfill\rule{1ex}{1ex}}
\def\0{\emptyset}

\def\q{\hfill\rule{1ex}{1ex}}

\usepackage{amsthm}

\newtheorem{theorem}{Theorem}[section]
\newtheorem{definition}[theorem]{Definition}
\newtheorem{lemma}[theorem]{Lemma}

\newtheorem{observation}[theorem]{Observation}
\newtheorem{cor}[theorem]{Corollary}
\newtheorem{prop}[theorem]{Proposition}

\newtheorem{conj}[theorem]{Conjecture}
\newtheorem{prob}[theorem]{Problem}

\usepackage{algorithm}
\usepackage{algpseudocode}
\usepackage{graphicx}
\usepackage{amsmath}
\usepackage{pstricks,pgf,subcaption,amssymb}
\begin{document}
\title{\bf Partite saturation number of cycles
	 }
\author[]{
Yiduo Xu\thanks{E-mail:\texttt{xyd23@mails.tsinghua.edu.cn}}\;}
\author[2]{
Zhen He\thanks{Corresponding author. E-mail:\texttt{zhenhe@bjtu.edu.cn}}}
\author[]{
Mei Lu\thanks{E-mail:\texttt{lumei@mail.tsinghua.edu.cn}}}

\affil[]{\small Department of Mathematical Sciences, Tsinghua University, Beijing 100084, China}
\affil[2]{\small School of Mathematics and Statistics, Beijing Jiaotong University, Beijing 100044, China.}
\date{}

\maketitle\baselineskip 16.3pt

\begin{abstract}
	A graph $H$ is said to be $F$-saturated relative to $G$, if $H$ does not contain any copy of $F$, but the addition of any edge $e$ in $E(G)\backslash E(H)$ would create a copy of $F$.  The minimum size of an $F$-saturated graph relative to $G$ is denoted by $sat(G,F)$. Let $K_k^n$ be the complete $k$-partite graph containing $n$ vertices in each part and $C_\ell$ be the cycle of length $\ell$. In this paper we give an asymptotically tight bound of $sat(K_k^n,C_\ell)$ for all $ \ell \geq 4, k \geq 2$ except $(\ell,k)=(4,4)$. Moreover, we determined the exact value of $sat(K_k^n,C_\ell)$ for $ k>\ell=4 $ and $5 \geq \ell>k \geq 3$ and $(\ell,k)=(6,2)$.	
	\end {abstract}
	
	
	
\section{Introduction}

	In this paper we only consider finite, simple and undirected graphs. For a graph $G$, we use $V(G)$ to denote the vertex set of $G$, $E(G)$ the edge set of $G$, $|G|$ the order of $G$ and $e(G)$ the size of $G$. For a positive integer $k$, let $[k]:=\{1,2,\ldots,k\}$. We say a graph $G$ is a $k$-partite graph, if $V(G)$ can be partition into $k$ parts $V_1,\ldots,V_k$ such that each of $V_i$ is an independent set. For a $k$-partite graph $G$, we simply write $G=V_1 \cup V_2 \cup \cdots \cup V_k$, where $V_i$ is the $i$-th part of $G$. Still we abbreviate $G$ as bipartite if $k=2$, tripartite if $k = 3$ and multipartite if $k \geq 3$. Denote by $C_r,P_r$ the $r$-vertices cycle and path, respectively. The length of a path (cycle) is the vertex number of the path (cycle). 
	
	 For a graph $G$ and $e \not \in E(G)$, we use $G+e$ to denote the graph obtained by adding the edge $e$ into $G$. Similarly, given a graph $G$ and $u \in V(G), e \in E(G)$, let $G - u$ (resp. $G - e$) be the graph obtained by deleting the vertex $u$ and all edges relating to $u$ (resp. deleting the edge $e$) in $G$.  For a vertex set (resp. an edge set) $A$, let $G[A]$ be the subgraph of $G$ induced by $A$. 
	
	For a graph $G$ and $u \in V(G)$, let $N_G(u)=\{ v : uv \in E(G) \}$ be the \textit{neighbourhood} of $u$, and $d_G(u)=|N_G(u)|$  the degree of $u$. We use $\delta(G)$  to denote the minimum  degree of  $G$. Without confusion, we
abbreviate as $N(u)$, $d(u)$ and $\delta$, respectively. For a vertex set $A $, let $N(A)= \bigcup \limits_{u \in A} N(u) \, \setminus \, A$. Let $N[u]=N(u) \cup \{u\}$ and $N[A] = N(A) \cup A$. If $G$ is a partite graph, let $N_i(u)=N(u) \cap V_i$ and $N_i(A)=N(A) \cap V_i$. The \textit{distance} $d(u,v)$ between two vertices $u,v$ is the number of edges contained in the shortest path connecting $u$ and $v$. If $u,v$ are not in the same connected components, then $d(u,v) = \infty$. Let $d(u,v)=0$ if $u=v$. For a given partite graph $G$, we say  a nonedge $uv$ is \textit{admissible} if $u,v$ lie in different parts of $G$. We use $diam_p(G)$ to denote the maximum distance between two vertices not lie in the same part of $G$. 
	
	Given graphs $G$ and $F$, we say a subgraph $H$ of $G$ is \textit{$F$-saturated relative to $G$}, if $H$ does not contain any copy of $F$ but the addition of any edge $e$ in $E(G)\backslash E(H)$ would create a copy of $F$. The \textit{saturation number} of $F$ relative to $G$ is denoted by
	\begin{equation*}
	   sat(G,F)=\min \{ e(H):H  \; \text{is}  \; F\text{-saturated relative to } G \} \, .
	\end{equation*}
Let $Sat(G,F) = \{H : H \; \text{is}  \; F\text{-saturated relative to } G \text{ and } e(H)= sat(G,F) \}$. If $G=K_n$, we abbreviate $sat(G,F)$ as $sat(n,F)$. The first saturation problem was studied in 1964 by Erd\H os, Hajnal and Moon \cite{EHM} who proved that $sat(n,K_r)=(r-2)(n-1)-\frac{(r-2)(r-1)}{2}$. For readers interested in saturation problem, we refer to the survey \cite{survey}.
	
	Finding cycle saturation number $sat(n,C_k)$ is an interesting problem in extremal graph theory. The saturation number $sat(n,C_3)$ is  given in \cite{EHM} mentioned above. Ollmann \cite{OLL} and Tuza \cite{TUZ} determined $sat(n,C_4)= \lfloor \frac{3n-5}{2} \rfloor$ for $n \geq 5$. Fisher, Fraughnaugh and Langley \cite{FI2} derived an upper bound of $sat(n,C_5)$ and was confirmed to be the exactly value of $sat(n,C_5)$ by Chen \cite{Chen1,Chen2} later, that is $sat(n, C_5)= \lceil \frac{10(n-1)}{7} \rceil$ for $n \geq 21$. Gould, Luczak,  Schmitt \cite{Gou} and Zhang, Luo,  Shigeno \cite{Zhang} proved that $ \lceil \frac{7n}{6} \rceil -2 \leq  sat(n,C_6) \leq \lfloor \frac{3n-3}{2} \rfloor$ for $n \geq 9$. Recently Lan,  Shi,  Wang and Zhang \cite{LAN} showed that $sat(n,C_6)=\frac{4}{3}n+O(1)$ as $n \geq 9$. For $k \geq 7$, F{\"u}redi and Kim \cite{Fur} showed that $\left( 1 + \frac{1}{k+2}\right) n -1 \leq sat(n,C_k) \leq \left( 1+ \frac{1}{k-4} \right) + \binom{k-4}{2}$ for $n \geq 2k-5$.
	
	Another interesting problem is the partite saturation number, that is to determine $sat(G,F)$ for different $F$ where $G$ is a partite graph. Let $K_{n_1,n_2,\ldots,n_k}$ be the complete $k$-partite graph containing $n_i$ vertices in the $i$-th part. If $n_1=n_2=\cdots=n_k=n$, we abbreviate it as $K_k^n$. Note that $sat(K_k^1,F)=sat(k,F)$, the study of partite saturation number would help us to determine general saturation numbers. 
	
	In this paper we focus on the partite saturation number of cycles $sat(K_k^n,C_{\ell})$ for all $\ell \geq 4$ and $k \geq 2$. When $k=2$ we only consider $\ell$ being an even integer since there is no odd cycle in a bipartite graph. An easy observation on a cycle-saturated partite graph $G$ is that $G$ must be connected. Otherwise, the addition of two vertices from different connected components and different parts of $G$ would not form a cycle, a contradiction. Hence we have $sat(K_k^n,C_{\ell}) \geq kn-1$. Moreover, $diam_p(G) \leq \ell-1$.
	
	We first focus on bipartite saturation problem. Bollob\'as \cite{BOL} and Wessel \cite{WES} determined that $sat(K_{(n,n)},K_{(s,t)})=n^2-(n-s+1)(n-t+1)$, where $K_{(a,b)}$ is the complete bipartite graph with $a$ vertices in the first part and $b$ vertices in the second part. Moshkovitz and Shapira \cite{MOS} conjectured that $sat(K_{n,n},K_{s,t}) \geq (s+t-2)n-\lfloor (\frac{s+t-2}{2})^2 \rfloor$ for sufficiently large $n$. This conjectured has been studied by Gan, Kor\'andi, Sudakov \cite{GAN} and Chakraborti, Chen,  Hasabnis \cite{CHA} separately. The results of \cite{BOL,WES} showed that $sat(K_{n,n},C_4)=2n-1$. In \cite{DU}, Dudek and Wojda have asked the problem that determining $sat(K_{n_1,n_2},C_{2t})$ for $t>2$. We first answer this problem for $t=3$ and give an upper bound for $t \geq 4$.
\vspace{0.1em}	
	
\begin{theorem}\label{T11} For $\ell \geq 3$ and $n_1,n_2 \geq \ell+2$,  $sat(K_{n_1,n_2},C_{2 \ell}) \leq n_1+n_2 + \ell^2 - 3 \ell +1$. Moreover, $sat(K_{n_1,n_2},C_{6}) = n_1+n_2 + 1$.
\end{theorem}
\vspace{0.1em}
	

	For results of saturation in multipartite graph, Sullivan and Wenger \cite{SUL} gave the bounds of $sat(K_{n_1,n_2,n_3},K_{\ell,m,p})$ for $\ell \geq m \geq p $ and $n_1,n_2,n_3$ sufficiently large. Moreover, they determined the exact value of $sat(K_{n_1,n_2,n_3},K_{\ell,\ell,p})$ when $p=\ell, \ell-1$. He and Lu \cite{HE} determined the exact value of $sat(K_{n_1,n_2,n_3},tK_{\ell,\ell,\ell})$ for $t,\ell \ge 1$ and $n_1,n_2, n_3$ sufficiently large. Ferrara, Jacobson, Pfender and  Wenger \cite{FJP} determined $sat(K_k^n,K_3)$ for $k \geq 3$ and $n \ge 100$. Roberts \cite{ROB} showed that
    $sat(K_4^n,K_4) = 18n - 21$ for sufficiently large $n$. Gir\~ ao, Kittipassorn and  Popielarz \cite{GIRAO} studied $sat(K_k^n,K_r) $ for $k \geq r \geq 3$ and sufficiently large $n $. Since \cite{FJP,GIRAO,SUL} have proved that $sat(K_k^n,C_3)=3(k-1)n-6$ for $k \geq 3$ and large $n$, we only consider $C_{\ell}$-saturated $k$-partite graph with $\ell \geq 4 $ and $ k \geq 3$. We give the general bounds of $sat(K_k^n,C_{\ell})$.

\vspace{0.1em}	
	
\begin{theorem}\label{T12} For any $ \ell \geq 4,  k \geq 3$ and $n \geq \frac{\ell}{k}$, $sat(K_k^n,C_{\ell}) \geq kn$.
\end{theorem}
\noindent{\bf Proof. }Since a cycle-saturated multipartite graph must be connected, we have $sat(K_k^n,C_{\ell}) \geq kn-1$. Let $G$ be a $C_{\ell}$-saturated $k$-partite graph with $e(G)=sat(K_k^n,C_{\ell})$. If $e(G)=kn-1$, then $G$ is a tree. If there exists $u \in V(G)$ such that $u$ has two neighbours $v,w$ lie in different parts, then  $G+vw$ does not contain $C_{\ell}$ as a subgraph. Thus the neighbours of any vertex $u \in V(G)$ must  lie in the same part which shows that $G$ is bipartite, a contradiction with $k \geq 3$. Hence  we have $e(G) \geq kn$. \qed
\vspace{0.8em}
\begin{theorem}\label{T13} For any  $\ell > k \geq 3$, $\ell \geq 6$ and $n \geq \lfloor \frac{\ell-2}{2} \rfloor$, 
\begin{equation*}
\begin{aligned}
sat(K_k^n,C_{\ell})\;  & \leq
\; kn- \ell+1 + \lfloor \frac{\ell-2}{2} \rfloor ^2  + 2 \left( \ell-1-2 \, \lfloor \frac{\ell-2}{2} \rfloor \right) \lfloor \frac{\ell-2}{2} 
\rfloor
\\
& =
\left\{ \;
\begin{aligned}
& kn + \frac{\ell^2-2\ell-11}{4}  \, , & \; & \ell  \text{ is odd} \, ; \\
& kn + \frac{\ell^2-4\ell}{4}   \, , & \; & \ell \text{ is even}  \, . 
\end{aligned}
\right. 
\end{aligned}
\end{equation*}
\end{theorem}
\vspace{0.1em}

\begin{theorem}\label{T14} For any  $k \geq \ell \geq 5$ and $n \geq \lceil \frac{k+\ell-2}{2} \rceil$, 
\begin{equation*}
\begin{aligned}
sat(K_k^n,C_{\ell})\;  & \leq
\; kn+k+ \ell-5 + \lfloor \frac{\ell-4}{2} \rfloor ^2  + 2 \left( \ell-4-2 \, \lfloor \frac{\ell-4}{2} \rfloor \right) \lfloor \frac{\ell-4}{2} 
\rfloor
\\
& =
\left\{ \;
\begin{aligned}
& kn+k +\frac{\ell^2-2\ell-15}{4} \, , & \; & \ell  \text{ is odd} \, ; \\
& kn+k +\frac{\ell^2-4 \ell -4}{4}  \, , & \; & \ell \text{ is even}  \, .
\end{aligned}
\right.
\end{aligned}
\end{equation*}
\end{theorem}
\vspace{0.1em}

	Theorems \ref{T13} and \ref{T14} give an upper bound of $sat(K_k^n,C_{\ell})$ for all $\ell \geq 4, k \geq 3$ except for $\ell=4$ and $\ell =5, 3 \leq k \leq 4$. The three results below show the exact value of $sat(K_k^n,C_{\ell})$ under cases $5 \geq \ell > k \geq 3$. In \cite{SUL} the authors gave an easy proof on $sat(K_3^n,C_{4})$. In this paper we give a different structural proof of $sat(K_3^n,C_{4})$, not only show the exact value of it but also show that the extremal graph of $sat(K_3^n,C_{4})$ is in $\Omega^{(4,3,n)}$, where $\Omega^{(4,3,n)}$ would be defined in Section 3.
\vspace{0.1em}

\begin{theorem}\label{T15} For any $n \geq 2$,  $$sat(K_3^n,C_{4}) = 3n \, , \quad Sat(K_{3 }^{ n}, C_4) = \Omega^{(4,3,n)}.$$
\end{theorem}
\vspace{0.1em}

\begin{theorem}\label{T16}   $sat(K_3^n, C_5)  = \left\{ \begin{matrix}
  \; 6 \, ,  \quad  \quad \; & \quad  n=2 \, ,  \\
  \; 3n+1 \, , & \quad n \geq 3 \, .
\end{matrix} \right. $
\end{theorem}
\vspace{0.0em}

\begin{theorem}\label{T17} For any $n \geq 10$,  $sat(K_4^n, C_5) =  4n+2 $.
\end{theorem}
\vspace{0.1em}
	
	Theorem \ref{T15} shows that the bounds in Theorem \ref{T12} are tight when $ \ell =4, k =3$. Theorems \ref{T11} - \ref{T17} immediately give Corollary \ref{C18}. 
\vspace{0.1em}
\begin{cor}\label{C18}
For any fixed $\ell \geq 5, k \geq 3$ or $\ell=4, k=3$ or even $ \ell \geq 6, k=2$, 
\begin{equation*}
 sat(K_k^n, C_{\ell}) = kn +  O(k+\ell^2) \, ,
\end{equation*}
as $n \rightarrow \infty$.
\end{cor}
\vspace{0.1em}
	
	For the rest case $k \geq \ell=4$, we have following results on $C_4$-saturated $k$-partite graph. When $n=1$, such value is equal to $sat(k,C_4)= \lfloor \frac{3k-5}{2} \rfloor$, answering that studying partite saturation number is useful to determine original saturation number.
	
\begin{theorem}\label{T19} For any $k \geq 5$ and $n \geq 1$, 
\begin{equation*}
sat(K_k^n, C_4) \; =\; \lfloor  \frac{3(k-1)n -2}{2} \rfloor = 
\left\{ \; 
\begin{aligned}
& \frac{3}{2}(k-1)n - \frac{3}{2} \, , & \; & k \text{ is even} \, , \, n \text{ is odd} \, ; \\
\\
& \frac{3}{2}(k-1)n - 1 \, , & \; &  \text{otherwise}  \, . \\
\end{aligned}
\right. 
\end{equation*}
Moreover, $ \lfloor \frac{9}{2}n-1 \rfloor \leq sat(K_4^n, C_4) \leq 5n-1$.
\end{theorem}
		
	For $k \geq 3$ and $\ell$ is much greater than $k$, we give a better construction in order to prove that $sat(K_k^n,C_{\ell}) = kn + O(\ell)$.
	\vspace{0.1em}
	
\begin{theorem}\label{T110} For any fixed $k \geq 4$, $\ell \geq 60k+12$ and $n$ sufficeintly large, $$ sat(K_k^n,C_{\ell}) \leq  k(n-1) + 6 \, \lceil \,  \frac{\ell}{5}  \, \rceil \,.$$
Moreover, $sat(K_3^n, C_{\ell}) \leq  3n + 13 \, \lceil \,  \frac{\ell}{5}  \, \rceil $ while $n,\ell$ satisfy the same condition.
\end{theorem}
\vspace{0.1em}

	The existence of $\lim \limits_{n \rightarrow \infty} \frac{sat(n,F)}{n}$ is the most important problem in saturation number. Tuza \cite{TUZA} proved that $sat(n,F)=O(n)$ for all $F$ and has conjectured that there exists a constant $c=c(F)$ such that $sat(n,F)=cn+o(n)$. This conjectured has been confirmed to be correct for many specific graphs, readers can see them in \cite{survey}. Such a conjecture can be asked to partite saturation version.
\vspace{0.1em}
	
\begin{conj}\label{C111} For any fixed $k \geq 2$ and $F$ being a subgraph of $K_k^n$, $\lim \limits_{n \rightarrow \infty} \frac{sat(K_k^n,F)}{kn}$ exists.
\end{conj}
 
	The partite saturation results introduced above and all Theorems given in this paper have confirmed that this conjecture is correct for specific graphs. Combining results of \cite{BOL,FJP,GIRAO,SUL,WES} we have following Corollary.
\vspace{0.1em}
	
\begin{cor}\label{C112} For any fixed $\ell , k \geq 3$ and even $\ell \geq 4,k=2$ except for $(\ell,k)=(4,4)$, $\lim \limits_{n \rightarrow \infty} \frac{sat(K_k^n,C_{\ell})}{kn}$ exists.
\end{cor}
\vspace{0.1em}

	From the table below, readers can visualize all Theorems and Corollaries in this section. The rest of this paper is organized as follows. In Section 2 we prove Theorem \ref{T11}.  In Section 3 we construct saturated subgraphs of $K_k^n$. We also prove Theorems \ref{T13}, \ref{T14} and \ref{T110} in this Section. In Section 4, we prove Theorems \ref{T15} and \ref{T19}. In Section 5, we prove Theorems \ref{T16} and \ref{T17}.  Finally, we conclude the paper in Section 6 with some open problems.
	
\begin{table}[H]
\renewcommand\arraystretch{1.8}
\centering
\large
\resizebox{\textwidth}{!}{
\begin{tabular}{|c|c|c|c|c|c|c|}
\hline
 \textcolor{red!60!black!100!}{$sat(K_k^n,C_\ell)$}    &   \textcolor{red!60!black!100!}{$k=2$} &   \textcolor{red!60!black!100!}{$k=3$}   &  \textcolor{red!60!black!100!}{\quad \quad $k=4$ \quad \quad}  & \multicolumn{3}{c|}{\textcolor{red!60!black!100!}{$k \geq 5$}}\\ \hline
 \textcolor{red!60!black!100!}{$\ell =3 $ } & \XSolid  &  \multicolumn{5}{c|}{ \; \quad $ =3(k-1)n-6 $ \textcolor{green!30!black!100!}{ \cite{FJP,GIRAO,SUL}} \quad \; } 
\\ \hline
 \textcolor{red!60!black!100!}{$\ell=4$} & $=2n-1 $ \textcolor{green!30!black!100!}{ \cite{BOL,WES} } & $=3n$ & \makecell[c]{ $\leq 5n-1$ \\ $\geq \lfloor \frac{9}{2}n-1 \rfloor$  } &  \multicolumn{3}{c|}{$=\lfloor \frac{3(k-1)n-2}{2} \rfloor$} \\ \hline
\textcolor{red!60!black!100!}{$\ell=5$} & \XSolid    & $=3n+1$ & $=4n+2$ &  \multicolumn{3}{c|}{$ \leq kn+k$} \\ \hline
\multirow{4}{*}{\textcolor{red!60!black!100!}{$  \ell \geq 6$}} & \multicolumn{1}{c|}{\textcolor{blue!60!black!100!}{$ \ell=6$}} & \multicolumn{4}{c|}{\textcolor{blue!60!black!100!}{$ \ell \geq k$}} &  \multicolumn{1}{c|}{\textcolor{blue!60!black!100!}{$ k > \ell$}} \\ \cline{2-7}
    &  \multicolumn{1}{c|}{ \;  $ =2n+1  $ \quad  } &  \multicolumn{4}{c|}{\; $< kn + \frac{\ell^2}{4}  $ \quad }  &   \multicolumn{1}{c|}{ \;  $ < k(n+1) + \frac{\ell^2}{4}   $ \quad  }  \\ \cline{2-7}
      & \multicolumn{1}{c|}{\textcolor{blue!60!black!100!}{$ \ell \geq 8$, even}} &   \multicolumn{4}{c|}{\textcolor{blue!60!black!100!}{$ \ell \geq 60k+11$} }    &  \, \\ \cline{2-6}
      & \multicolumn{1}{c|}{ \;  $ < 2n+\frac{\ell^2}{4}  $  \quad } &  \;  $ \leq 3n + 13 \, \lceil \,  \frac{\ell}{5}  \, \rceil $ \quad  & \multicolumn{3}{c|}{ \; $\leq  k(n-1) + 6 \, \left\lceil \,  \frac{\ell}{5}  \, \right\rceil $ \quad }  &   \,  \\
\hline 

\end{tabular}

}

\caption{\centering The values or bounds of $sat(K_k^n,C_\ell)$ for all $\ell\geq 4, k \geq 2$.} 
\end{table}	
	
\section{Bipartite saturation number of even cycles}	
We begin this section by an easy observation, and we would construct $C_{2\ell}$-saturated bipartite graph for $\ell \geq 3$.

\begin{observation}\label{O21} For $\ell \geq 3$, let $H=K_{\ell,\ell-1}= V_1 \cup V_2$ where $|V_1|=\ell,|V_2|=\ell-1$. Then for any different $u,v \in V(H)$, \\
\indent (1) If $u,v \in V_1$, then there exists $P_{2\ell-1},P_{2\ell-3}$ connecting $u,v$ in $H$. \\
\indent (2) If $u \in V_1, v \in V_2$, then there exists $P_{2\ell-2}$ connecting $u,v$ in $H$.
\end{observation}
\vspace{0.0em}

\begin{definition}\label{D22}
Let $\ell\geq 3$ and $n_1,n_2 \geq \ell+2$. Let $\mathscr{G}_{n_1,n_2}^{\ell} = V_1 \cup V_2$ be a bipartite graph such that $|V_1|=n_1,|V_2|=n_2$. Let $A_i \subseteq V_i $ be a $(n_i-\ell)$-element vertex set and $B_i = V_i \setminus A_i$. Pick distinct vertices $x_1,x_2 \in B_1$ and $y_1 \in B_2$. Then the edge $e$ in $G$ can only be\,:\\
\indent (1) $e=uy_1$ for any $u \in A_1$\,; \\
\indent (2) $e=vx_1$ for any $v \in A_2$\,; \\
\indent (3) $e=y_1x_2$\,; \\
\indent (4) $e=xy$ for any $u \in B_1, v \in B_2 \setminus\{ y_1\}$. 
\end{definition}

\begin{figure}[H]
\centering
\begin{tikzpicture}[scale=.45]

\draw (-9.8,0) arc(0:360:1.2cm and 3.5cm) ;
\draw (-5.8,0) arc(0:360:1.2cm and 3.5cm) ;

\draw (-9.8,-8) arc(0:360:1.2cm and 3cm) ;
\draw (-5.8,-8.4) arc(0:360:1.2cm and 2.6cm) ;

\draw (-7,0) node[align=center]{$A_2$};
\draw (-11,0) node[align=center]{$A_1$};

\draw (-7,-8.4) node[align=center]{$B_2'$};
\draw (-11,-9) node[align=center]{$B_1$};

\filldraw (-10.99,-7) circle (5pt);
\filldraw (-6.99,-5) circle (5pt);
\filldraw (-10.99,-6) circle (5pt);

\draw (-11.4,-7.6) node[align=center]{$x_2$};
\draw (-6.2,-5.5) node[align=center]{$y_1$};
\draw (-11.4,-6.4) node[align=center]{$x_1$};

\draw[thick] (-10.99,-7)--(-6.99,-5) (-9.8,-8) -- (-8.2,-8)(-9.8,-8.02) -- (-8.2,-8.02) (-9.8,-7.98) -- (-8.2,-7.98);

\draw[thick] (-10.45,-3) --(-6.99,-5) (-7.55,-3) --(-10.99,-6);

\end{tikzpicture}\\

\caption{\centering $\mathscr{G}_{n_1,n_2}^{\ell} = V_1 \cup V_2$. $V_i = B_i \cup A_i$, $B_2'=B_2 \setminus \{y_1\}$ such that $|B_1|=\ell$, $|B_2'|=\ell-1$. The solid line represents the complete connection between vertices. }

\end{figure}
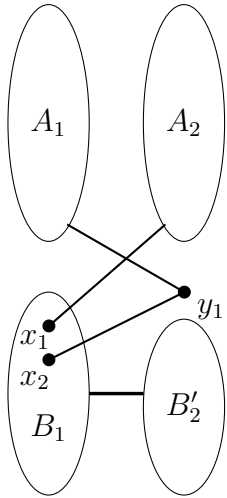

	Readers can visualize $\mathscr{G}_{n_1,n_2}^{\ell}$ in Figure 1. It is clear that $e(\mathscr{G}_{n_1,n_2}^{\ell})=n_1+n_2+\ell^2-3\ell+1$.
\vspace{0.0em}

\begin{prop}\label{P23} 
$\mathscr{G}_{n_1,n_2}^{\ell}$ is $C_{2 \ell}$-saturated.
\end{prop}
\noindent{\bf Proof. }Obviously all vertices in $A_1 \cup A_2 \cup \{y_1\}$ are not in cycles, so $\mathscr{G}_{n_1,n_2}^{\ell}$ is $C_{2\ell}$-free. To prove that $\mathscr{G}_{n_1,n_2}^{\ell}$ is $C_{2 \ell}$-saturated we only need to show that for any admissible nonedge $uv$ there exists a $P_{2\ell}$ connecting $u,v$. Let $H=\mathscr{G}_{n_1,n_2}^{\ell}[B_1 \cup B_2'] = K_{\ell,\ell-1}$.

	Let $uv$ be an admissible nonedge. If $u \in A_1, v \in A_2$, by Observation \ref{O21} there exists $P_{2\ell-3}$ connecting $x_1,x_2$ in $H$, together with $u,v,y_1$ there exists $P_{2\ell}$ connecting $u,v$. For the rest cases (1) $u \in A_1, v \in B_2'$, (2) $u \in B_1 \setminus\{x_1\}, v \in A_2$ and (3) $u \in B_1 \setminus\{x_2\}, v =y_2$ we are also done by similar arguments using Observation \ref{O21}. Hence $\mathscr{G}_{n_1,n_2}^{\ell}$ is $C_{2 \ell}$-saturated. \qed
	
\vspace{0.3em}

\begin{prop}\label{P24} 
For $n_1,n_2 \geq 5$, $sat(K_{n_1,n_2},C_6) \geq n_1+n_2 +1$.
\end{prop}
\noindent{\bf Proof. }Let $G=V_1 \cup V_2$ be a $C_6$-saturated bipartite graph with $|V_i|=n_i$. If $e(G)=n_1+n_2$, since $G$ is connected $G$ contains exactly one even cycle $C_{2r}=x_1y_1 \ldots x_ry_r$ where $x_i \in V_1, y_i \in V_2$. Since $G$ is $C_6$-saturated we have $r \neq 3$.

	If $r \geq 5$, then the addition of $x_1y_2$ in $G$ would only create $C_4$ and cycles of length at least 8, a contradiction. If $r=2$ or $r=4$, since $n_1 \geq 5$ there must exists $x_i \neq x^* \in V_1$ for $i \in [r]$ such that $|N(x^*) \cap \{y_1,\ldots y_r\}| =1$. Assume that $N(x^*) \cap \{y_1,\ldots y_r\} = \{y_1\}$, then the addition of $x^*y_2$ would only create $C_4$ or $C_8$, a contradiction. Hence $e(G) \geq n_1+n_2 +1$. \qed
\vspace{1em}
	
	\noindent{\bf Proof of Theorem \ref{T11}. }By Propositions \ref{P23} and \ref{P24} we are done. \qed

\section{Constructions of saturated subgraphs of $K_{k}^{ n}$ }	

In this section, we construct several graph families to show the upper bounds of our saturation results. We begin from a simple proposition.

\begin{prop}\label{P31} Let $G=V_1 \cup \cdots \cup V_k$ be a $k$-partite graph. If there exist $A \subseteq V(G)$ and $x_1,\ldots,x_k \in A$ such that:\\
\indent (1) $N(u_i)=\{x_i\}$ for any $u_i \in V_i \setminus A$ and $i\in[k]$, and there exists a $P_{\ell-1}$ connecting $x$ and $x_i$ for any $x \in A\setminus (V_i\cup\{x_i\})$; \\
\indent (2) There exists a $P_{\ell-2}$ connecting $x_i$ and $x_j$ for any $i,j\in [k]$; \\
\indent (3) $G[A]$ is $C_{\ell}$-saturated.\\
\noindent Then $G$ is  $C_{\ell}$-saturated.
\end{prop}
\vspace{0.0em}

\noindent{\bf Proof. }Obviously, $G$ is $C_{\ell}$-free. By (3), we only need to consider an admissible nonedge who has at least one endpoint in $V(G) \setminus A$. Let $u_iw$ be an admissible nonedge with $u_i \in V_i \setminus A$, where $i\in [k]$.
 If $w\in V_j\setminus A$ for some $j\in [k]\setminus\{i\}$, by (1) and (2), there exists a $P_{\ell}$ connecting $u_i$ and $w$. If $w \in A$, then $w\notin V_i$. By (1), there exists a  $P_{\ell}$ connecting $u_i$ and $w$. Hence  $G+u_iw$ would create a $C_{\ell}$ in each case. Thus $G$ is  $C_{\ell}$-saturated. \qed
\vspace{0.8em}

	Let $G=V_1 \cup \cdots \cup V_k$ be a $k$-partite graph and $x_1,\ldots,x_k\in V(G)$. We call $(x_1,\ldots,x_k)$  a \textit{ good pair} of $G$ if they do not lie in the same part of $G$.

\subsection{Construction I}

We first construct a graph family $\gamma^{(\ell,k)}$. 
\vspace{0.2em}

\begin{definition}\label{D32}
Let $\ell > k \geq 3$. \\
\indent (1) If $ \ell \leq 5$ and $(\ell, k ) \neq (5,3)$, let $\gamma^{(4,3)}=K_3$, $\gamma^{(5,4)}=K_4$. \\
\indent (2) If $\ell \geq 6$, let
\begin{equation*}
\gamma^{(\ell,k)} = K_{\lfloor \frac{\ell-2}{2} \rfloor \, , \, \lfloor \frac{\ell-2}{2} \rfloor\,  ,  \, \left( \ell-1 - 2 \lfloor \frac{\ell-2}{2} \rfloor \right) } = \left\{
\begin{aligned}
& K_{ \frac{\ell-2}{2}  ,  \frac{\ell-2}{2} ,\, 1 } \, , & \; & \ell \text{ is even} \, ,  \\
& K_{ \frac{\ell-3}{2} , \frac{\ell-3}{2}  ,\, 2 } \, , & \; & \ell \text{ is odd}    \, .
\end{aligned}
\right.
\end{equation*}
\end{definition}

\vspace{0.0em}

	Not difficult to see that any pairs of vertices $u,v$ in $\gamma^{(\ell,k)}$ are connected by $P_{\ell-1}$ and $P_{\ell-2}$ for $(\ell,k) \neq (5,3)$ and $\gamma^{(\ell,k)}$ is $C_{\ell}$-saturated because $|\gamma^{(\ell,k)}| = \ell-1$. Since $\ell-1 \geq k$, the good pair of $\gamma^{(\ell,k)}$ exists. The next lemma could help us to construct $\Omega^{(\ell,k,n)}$.
\vspace{0.0em}

\begin{lemma}\label{L33}
For $ k \geq 3$, let $G$ be a $k$-partite graph and $(x_1,\ldots,x_k)$ be a
 good pair of $G$.  Then there exists a bijection $f \,: \, \{x_1,\ldots,x_k \} \rightarrow [k] $ such that $f(x_i) \neq \pi(x_i)$, where $\pi(x_i)$ is the part of $G$ containing $x_i$.
\end{lemma}

\noindent{\bf Proof. }Assume that the number of parts containing $x_1,\ldots,x_k$ is $r$, then $r\ge 2$. Assume, without loss of generality, that $x_i \in V_i$ for $i \in [r]$, and $x_j \in V_1 \cup \cdots \cup V_r$ for $j \in [k] \setminus [r]$. Denote
\begin{equation*}
\begin{aligned}
& f(x_1)=r \, , \\
& f(x_i)=i-1 \, , \;  \, i=2,\ldots,r \, , \\
& f(x_j)=j \, , \; 
 \, j=r+1,\ldots,k \, .
\end{aligned}
\end{equation*}
It is easy to check that $f$ satisfies the property. \qed
\vspace{0.3em}

\begin{definition}\label{D34}
For $n \geq  \lfloor \frac{\ell-2}{2} \rfloor$, $\ell > k \geq 3$ and $(\ell,k) \neq (5,3)$, let $W^{(\ell,k,n)}$ be a graph   satisfies the following properties:  \\
\indent \text{(i)} $W^{(\ell,k,n)}$ is a $k$-partite graph such that each part of $W^{(\ell,k,n)}$ has exactly $n$ vertices. \\
\indent \text{(ii)} $V(W^{(\ell,k,n)})=A \cup B$ with $A \cap B = \emptyset$, $W^{(\ell,k,n)}[A]=\gamma^{(\ell,k)}$  and $W^{(\ell,k,n)}[B]$ is an empty graph. \\
\indent \text{(iii)} For any vertex $v \in B$, $d_{W^{(\ell,k,n)}}(v)=1$, $N(v) \subseteq  \{x_1,\ldots,x_k \}$, where $(x_1,\ldots,x_k)$ is a good pair of $W^{(\ell,k,n)}[A]$. For different $u,v \in B$, if $u,v$ lie in the same part of $W^{(\ell,k,n)}$, then $N(u)=N(v)$; otherwise  $N(u) \cap N(v) = \emptyset$.
\end{definition}
\vspace{0.0em}

	Definition \ref{D34} (iii) can be written as a bijection  $f : \{x_1,\ldots,x_k \} \rightarrow 2^B$ satisfying $  \bigcup \limits_{i=1}^k f(x_i)= B$, $f(x_i) = V_{m_i} \cap B$ for some $m_i$, and $f(x_i) \cap f(x_j) = \emptyset $ for all $i \neq j$. Since $W^{(\ell,k,n)}$ is a $k$-partite graph, Lemma \ref{L33} guarantees that such $f$ exists. If we want to emphasize the correlation between $W^{(\ell,k,n)}$ and $f$, we can rewrite $W^{(\ell,k,n)}$ as $W^{(\ell,k,n)}_\pi$, where $\pi = \pi(f , (x_1,\ldots,x_k), \ell , k)$ is an index depending on $f, (x_1,\ldots,x_k), \ell , k $.
Let $\Omega^{(\ell,k,n)}= \{ W^{(\ell,k,n)}_\pi \}$. If we ignore the repeated isomorphic graphs in $\Omega^{(\ell,k,n)}$, then we have $\Omega^{(4,3,n)}=\{W^{(4,3,n)}\}$ (Figure 2) and $\Omega^{(5,4,n)}=\{W_{\pi_1}^{(5,4,n)} , W_{\pi_2}^{(5,4,n)}\}$ (Figure 3). Still we give an example of $W^{(6,5,n)}_\pi \in \Omega^{(6,5,n)}$ (Figure 4).
	
	For $(\ell,k)=(5,3)$ and $n \geq 2$, let $\Omega^{(5,3,n)} = \{ W_*^{(5,3,n)} \} $ where $W_*^{(5,3,n)}$ is defined in Figure 5. The rest of this subsection is the proof of Theorem \ref{T13}.
	\vspace{1em}
	
	\noindent{\bf Proof of Theorem \ref{T13}. }For $\ell > k \geq 3 , \ell \geq 6$ and  $n \geq \lfloor \frac{\ell-2}{2} \rfloor$, let $G \in \Omega^{(\ell,k,n)}$. It is easy to check that $G$ satisfies the condition in Proposition \ref{P31}. Thus $G$ is a $C_{\ell}$-saturated $k$-partite graph which implies $sat(K_k^n, C_{\ell}) \leq e(G) = kn- \ell+1 + \lfloor \frac{\ell-2}{2} \rfloor ^2  + 2 \left( \ell-1-2  \lfloor \frac{\ell-2}{2} \rfloor \right) \lfloor \frac{\ell-2}{2} \rfloor$.   \qed  
	\vspace{1em}

\begin{figure}[H]
\centering
\begin{tikzpicture}[scale=.44]

\draw (-10,0) arc(0:360:1cm and 3.5cm) ;
\draw (-6,0) arc(0:360:1cm and 3.5cm) ;
\draw (-2,0) arc(0:360:1cm and 3.5cm) ;

\draw (-3,0) node[align=center]{$B_3$};
\draw (-7,0) node[align=center]{$B_2$};
\draw (-11,0) node[align=center]{$B_1$};

\filldraw (-2.99,-5) circle (5pt);
\filldraw (-6.99,-5) circle (5pt);
\filldraw (-10.99,-5) circle (5pt);

\draw (-3.2,-4.4) node[align=center]{$x_3$};
\draw (-6.6,-4.4) node[align=center]{$x_2$};
\draw (-10.6,-4.4) node[align=center]{$x_1$};

\draw[dotted,thick] (-1, -5) arc(0:360:6cm and 1.2cm) ;

\draw (-2.99,-5)--(-6.99,-5) (-10.99,-5)--(-6.99,-5)  (-2.99,-5) arc(0:-180: 4cm and 0.7cm) (-10.99,-5);

\draw (-10.5,-3) --(-6.99,-5) (-6.5,-3) --(-2.99,-5) (-3.5,-3) --(-10.99,-5);

\end{tikzpicture}\\

\caption{\centering $W^{(4,3,n)} = V_1 \cup V_2 \cup V_3$, $V_i = B_i \cup \{x_i \}$, $A=\{x_1,x_2,x_3\}$. The solid line represents the complete connection between vertices, and the dotted ellipse represents $W^{(4,3,n)}[A]= \gamma^{(4,3)}=K_3$. }

\end{figure}
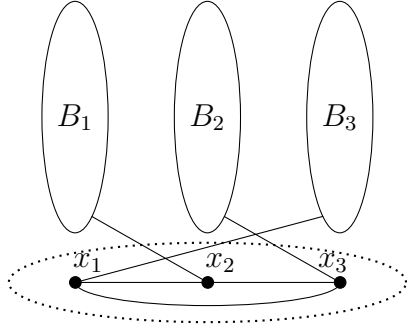

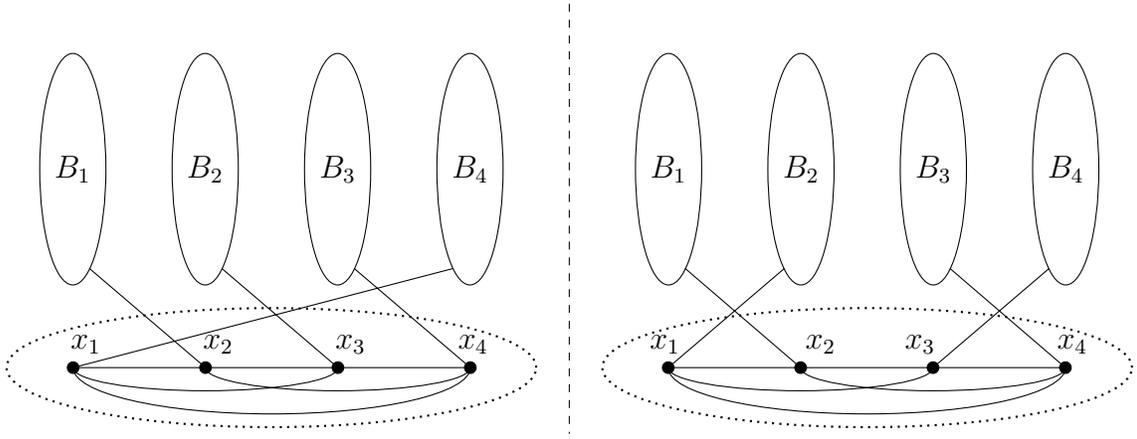
\begin{figure}[H]
\centering
\begin{tikzpicture}[scale=.44]

\draw (-14,1) arc(0:360:1cm and 3.5cm) ;
\draw (-10,1) arc(0:360:1cm and 3.5cm) ;
\draw (-6,1) arc(0:360:1cm and 3.5cm) ;
\draw (-2,1) arc(0:360:1cm and 3.5cm) ;

\draw (-3,1) node[align=center]{$B_4$};
\draw (-7,1) node[align=center]{$B_3$};
\draw (-11,1) node[align=center]{$B_2$};
\draw (-15,1) node[align=center]{$B_1$};

\filldraw (-2.99,-5) circle (5pt);
\filldraw (-6.99,-5) circle (5pt);
\filldraw (-10.99,-5) circle (5pt);
\filldraw (-14.99,-5) circle (5pt);

\draw (-2.9,-4.3) node[align=center]{$x_4$};
\draw (-6.6,-4.3) node[align=center]{$x_3$};
\draw (-10.6,-4.3) node[align=center]{$x_2$};
\draw (-14.6,-4.3) node[align=center]{$x_1$};

\draw[dotted,thick] (-1, -5) arc(0:360:8cm and 1.8cm) ;

\draw (-2.99,-5)--(-6.99,-5) (-10.99,-5)--(-6.99,-5) (-10.99,-5)--(-14.99,-5)  (-2.99,-5) arc(0:-180: 4cm and 0.7cm) (-10.99,-5) (-6.99,-5) arc(0:-180: 4cm and 0.7cm) (-14.99,-5)  (-2.99,-5) arc(0:-180: 6cm and 1.4cm) (-14.99,-5) ;

\draw (-10.5,-2) --(-6.99,-5) (-6.5,-2) --(-2.99,-5) (-3.5,-2) --(-14.99,-5) (-14.5,-2) --(-10.99,-5);

\draw[dashed] (0,6)--(0,-7);

\draw (16,1) arc(0:360:1cm and 3.5cm) ;
\draw (12,1) arc(0:360:1cm and 3.5cm) ;
\draw (8,1) arc(0:360:1cm and 3.5cm) ;
\draw (4,1) arc(0:360:1cm and 3.5cm) ;

\draw (15,1) node[align=center]{$B_4$};
\draw (11,1) node[align=center]{$B_3$};
\draw (7,1) node[align=center]{$B_2$};
\draw (3,1) node[align=center]{$B_1$};

\filldraw (2.99,-5) circle (5pt);
\filldraw (6.99,-5) circle (5pt);
\filldraw (10.99,-5) circle (5pt);
\filldraw (14.99,-5) circle (5pt);

\draw (15.2,-4.3) node[align=center]{$x_4$};
\draw (10.6,-4.3) node[align=center]{$x_3$};
\draw (7.6,-4.3) node[align=center]{$x_2$};
\draw (2.9,-4.3) node[align=center]{$x_1$};

\draw[dotted,thick] (17, -5) arc(0:360:8cm and 1.8cm) ;

\draw (2.99,-5)--(6.99,-5) (10.99,-5)--(6.99,-5) (10.99,-5)--(14.99,-5)  (10.99,-5) arc(0:-180: 4cm and 0.7cm) (2.99,-5)  (14.99,-5) arc(0:-180: 4cm and 0.7cm)(6.99,-5)  (14.99,-5)  (14.99,-5) arc(0:-180: 6cm and 1.4cm)  (2.99,-5);

\draw (6.5,-2) --(2.99,-5) (3.5,-2) --(6.99,-5) (14.5,-2) --(10.99,-5) (11.5,-2) --(14.99,-5);

\end{tikzpicture}\\

\caption{\centering $W_{\pi_1}^{(5,4,n)}$ (left) , $W_{\pi_2}^{(5,4,n)}$ (right). $W_{\pi_i}^{(5,4,n)} = V_1 \cup V_2 \cup V_3 \cup V_4$, $V_i = B_i \cup \{x_i \}$, $A=\{x_1,x_2,x_3,x_4\}$. The solid line represents the complete connection between vertices, and the dotted ellipse represents $W_{\pi_i}^{(5,4,n)}[A]= \gamma^{(5,4)}=K_4$.  }

\end{figure}

\begin{figure}[H]
\centering
\begin{tikzpicture}[scale=.44]

\draw (-14,1) arc(0:360:1cm and 3.5cm) ;
\draw (-10,1) arc(0:360:1cm and 3.5cm) ;
\draw (-6,0) arc(0:360:1cm and 4.5cm) ;
\draw (-2,-1) arc(0:360:1cm and 5.5cm) ;
\draw (2,-1) arc(0:360:1cm and 5.5cm) ;

\draw (1,-1) node[align=center]{$B_5$};
\draw (-3,-1) node[align=center]{$B_4$};
\draw (-7,0) node[align=center]{$B_3$};
\draw (-11,1) node[align=center]{$B_2$};
\draw (-15,1) node[align=center]{$B_1$};

\filldraw (-14.99,-4) circle (5pt);
\filldraw (-14.99,-6) circle (5pt);
\filldraw (-10.99,-4) circle (5pt);
\filldraw (-10.99,-6) circle (5pt);
\filldraw (-6.99,-6) circle (5pt);

\draw (-15.4,-3.3) node[align=center]{$x_1$};
\draw (-15.4,-6.7) node[align=center]{$x_2$};
\draw (-10.6,-3.3) node[align=center]{$x_3$};
\draw (-10.6,-6.7) node[align=center]{$x_4$};
\draw (-6.9,-6.9) node[align=center]{$x_5$};

\draw (-10.99,-4) -- (-14.99,-4) (-10.99,-6) -- (-14.99,-4) (-10.99,-4) -- (-14.99,-6) (-10.99,-6) -- (-14.99,-6) ;
\draw (-14.99,-4) -- (-6.99,-6) (-10.99,-4) -- (-6.99,-6) (-10.99,-6) -- (-6.99,-6)  (-6.99,-6) arc(0:-180: 4cm and 1.2cm) (-14.99,-6);

\draw[dotted,thick] (-15.99,-2.7) -- (-9.99,-2.7) -- (-5.99,-5.5) -- (-5.99,-7.5) -- (-15.99,-7.5) -- (-15.99,-2.7);

\draw (-14.99,-4) -- (-11.5,-2) (-10.99,-4) -- (-14.5,-2) (-14.99,-6)--(-7.5,-4) (-10.99,-6)  --(-3.8,-4.5) (0.6,-6) arc(0:-180: 3.8cm and 1.2cm) (-6.99,-6);

\end{tikzpicture}\\

\caption{\centering $W_{\pi}^{(6,5,n)} = V_1 \cup V_2 \cup V_3 \cup V_4 \cup V_5$, $V_1 = B_1 \cup \{x_1,x_2 \}$, $V_2 = B_2 \cup \{x_3,x_4 \}$, $V_3 = B_3 \cup \{x_5 \}$, $V_4=B_4$, $V_5=B_5$, $A=\{x_1,\ldots,x_5\}$. The solid line represents the complete connection between vertices, and the dotted pentagon represents $W_{\pi}^{(6,5,n)}[A]= \gamma^{(6,5)}=K_{2,2,1}$.   }

\end{figure}
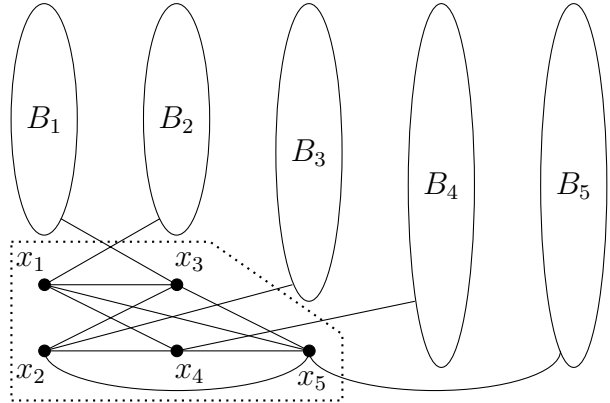

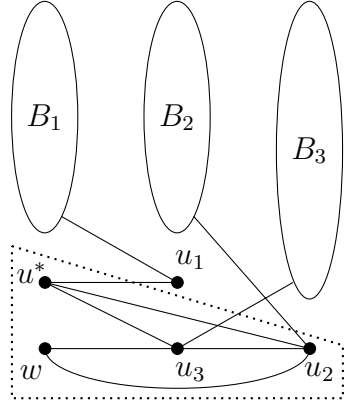
\begin{figure}[H]
\centering

\begin{tikzpicture}[scale=.44]

\draw (-14,1) arc(0:360:1cm and 3.5cm) ;
\draw (-10,1) arc(0:360:1cm and 3.5cm) ;
\draw (-6,0) arc(0:360:1cm and 4.5cm) ;

\draw (-7,0) node[align=center]{$B_3$};
\draw (-11,1) node[align=center]{$B_2$};
\draw (-15,1) node[align=center]{$B_1$};

\filldraw (-14.99,-4) circle (5pt);
\filldraw (-14.99,-6) circle (5pt);
\filldraw (-10.99,-4) circle (5pt);
\filldraw (-10.99,-6) circle (5pt);
\filldraw (-6.99,-6) circle (5pt);

\draw (-15.4,-3.7) node[align=center]{$u^*$};
\draw (-15.4,-6.7) node[align=center]{$w$};
\draw (-10.6,-3.3) node[align=center]{$u_1$};
\draw (-10.6,-6.7) node[align=center]{$u_3$};
\draw (-6.7,-6.7) node[align=center]{$u_2$};

\draw (-10.99,-4) -- (-14.99,-4) (-10.99,-6) -- (-14.99,-4)  (-10.99,-6) -- (-14.99,-6) ;
\draw (-14.99,-4) -- (-6.99,-6) (-10.99,-6) -- (-6.99,-6)  (-6.99,-6) arc(0:-180: 4cm and 1.2cm) (-14.99,-6);

\draw[dotted,thick] (-15.99,-2.9) -- (-5.99,-5.9) -- (-5.99,-7.5) -- (-15.99,-7.5) -- (-15.99,-2.9);

\draw (-6.99,-6) -- (-10.5,-2) (-10.99,-4) -- (-14.5,-2) (-10.99,-6)--(-7.5,-4) ;

\end{tikzpicture}\\
\caption{\centering $W_{*}^{(5,3,n)} = V_1 \cup V_2 \cup V_3 $, $|V_i|=n$. $V_1 = B_1 \cup \{u^*,w \}$, $V_2 = B_2 \cup \{u_1,u_3\}$, $V_3 = B_3 \cup \{u_2 \}$. The solid line represents the complete connection between vertices. }
\end{figure}

\subsection{Construction II}

In this subsection we first construct a graph family $\zeta^{(\ell,k)}$.

\begin{definition}\label{D35}
For $k \geq \ell \geq 5$, let $\zeta^{(\ell,k)} \subseteq V_1 \cup ... \cup V_5$ be a 5-partite graph on $k+\ell-2$ vertices. $V(\zeta^{(\ell,k)})=A \cup B \cup C$ such that $|A|=k, |B|=2$ and $|C|=\ell-4$. $A \subseteq V_1 \cup V_2$ and $\zeta^{(\ell,k)}[A]$  is an empty graph with half vertices in each part. $B=\{b_4,b_5\}$ such that $b_i \in V_i$ and $b_4b_5$ is an edge. $C \subseteq V_1 \cup V_2 \cup V_3 $ and 

\begin{equation*}
\zeta^{(\ell,k)}[C] = K_{\lfloor \frac{\ell-4}{2} \rfloor \, , \, \lfloor \frac{\ell-4}{2} \rfloor\,  ,  \, \left( \ell-4 - 2 \lfloor \frac{\ell-4}{2} \rfloor \right) } = \left\{
\begin{aligned}
& K_{ \frac{\ell-4}{2}  ,  \frac{\ell-4}{2} ,\, 0 } \, , & \; & \ell \text{ is even} \, ,  \\
& K_{ \frac{\ell-5}{2} , \frac{\ell-5}{2}  ,\, 1 } \, , & \; & \ell \text{ is odd}    \, .
\end{aligned}
\right.
\end{equation*}
That is, $\zeta^{(\ell,k)}[C]$ is a complete bipartite or tripartite graph (or is an isolated vertex) in the first three parts of $\zeta^{(\ell,k)}$. The rest of edges in $\zeta^{(\ell,k)}$ are all edges between $B$ and $A \cup C$.
\end{definition}
 	
 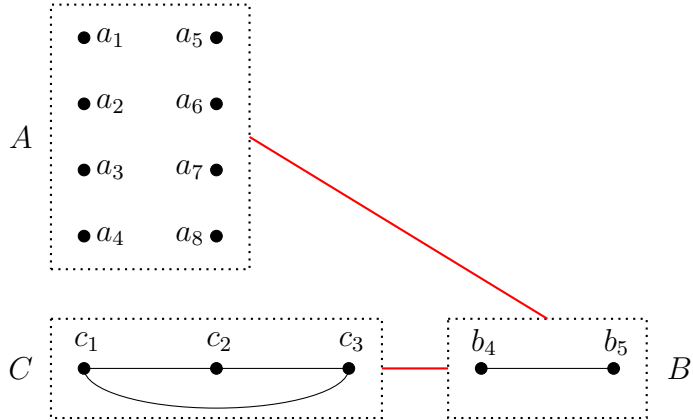
\begin{figure}[H]
\centering
\begin{tikzpicture}[scale=.44]



\filldraw (-14.99,-6) circle (5pt);
\filldraw (-10.99,-6) circle (5pt);
\filldraw (-6.99,-6) circle (5pt);
\filldraw (-2.99,-6) circle (5pt);
\filldraw (1.01,-6) circle (5pt);

\filldraw (-14.99,-2) circle (5pt);
\filldraw (-14.99,0) circle (5pt);
\filldraw (-14.99,2) circle (5pt);
\filldraw (-14.99,4) circle (5pt);

\filldraw (-10.99,-2) circle (5pt);
\filldraw (-10.99,0) circle (5pt);
\filldraw (-10.99,2) circle (5pt);
\filldraw (-10.99,4) circle (5pt);

\draw (3,-6) node[align=center]{$B$};
\draw (-16.9,-6) node[align=center]{$C$};
\draw (-16.9,1) node[align=center]{$A$};
\draw (-6.9,-5.2) node[align=center]{$c_3$};
\draw (-10.9,-5.2) node[align=center]{$c_2$};
\draw (-14.9,-5.2) node[align=center]{$c_1$};
\draw (-2.9,-5.2) node[align=center]{$b_4$};
\draw (1.1,-5.2) node[align=center]{$b_5$};

\draw (-14.19,-2) node[align=center]{$a_4$};
\draw (-14.19,0) node[align=center]{$a_3$};
\draw (-14.19,2) node[align=center]{$a_2$};
\draw (-14.19,4) node[align=center]{$a_1$};
\draw (-11.79,-2) node[align=center]{$a_8$};
\draw (-11.79,0) node[align=center]{$a_7$};
\draw (-11.79,2) node[align=center]{$a_6$};
\draw (-11.79,4) node[align=center]{$a_5$};

\draw (-10.99,-6) -- (-14.99,-6) (-2.99,-6) -- (1.01,-6);
\draw  (-10.99,-6) -- (-6.99,-6)  (-6.99,-6) arc(0:-180: 4cm and 1.2cm) (-14.99,-6);

\draw[dotted,thick] (-15.99,-4.5) -- (-5.99,-4.5) -- (-5.99,-7.5) -- (-15.99,-7.5) -- (-15.99,-4.5);

\draw[dotted,thick] (-3.99,-4.5) -- (2.01,-4.5) -- (2.01,-7.5) -- (-3.99,-7.5) -- (-3.99,-4.5);

\draw[dotted,thick] (-15.99,-3) -- (-9.99,-3) -- (-9.99,5) -- (-15.99,5) -- (-15.99,-3);

\draw[red,thick] (-5.99,-6) -- (-3.99,-6) (-0.99,-4.5) -- (-9.99,1) ;

\end{tikzpicture}\\

\caption{\centering $\zeta^{(7,8)} = V_1 \cup \cdots \cup V_5$, $V_1 = \{a_1,a_2,a_3,a_4,c_1\}$, $V_2 = \{a_5,a_6,a_7,a_8,c_2\}$, $V_3 =\{c_3 \}$, $V_4=\{b_4\}$, $V_5=\{b_5\}$, $A=\{a_1,..,a_8\}$, $B=\{b_4,b_5\}$, $C=\{c_1,c_2,c_3\}$. The red solid line represents the complete connection between vertices, and the black solid line represents edges.    }

\end{figure}

 	In Figure 6 we give an example of $\zeta^{(\ell,k)}$ for $(\ell,k)=(7,8)$. We now establish some properties on  $\zeta^{(\ell,k)}$.

\begin{prop}\label{P36} For $\ell \geq 6$, let $u,v$ be distinct vertices in $\zeta^{(\ell,k)}$, then we have\,:\\
\indent (1) If $u , v \in B \cup C$, then there exists  $P_{\ell-2}$, $P_{\ell-3}$ and $P_{\ell-4}$ connecting $u,v$ in $\zeta^{(\ell,k)}[B \cup C]$. \\
\indent (2) If $u \in A$, $v \in B$, then there exists  $P_{\ell-1}$ and $P_{\ell-2}$ connecting $u,v$. \\
\indent (3) If $u\in A$, $v \in A \cup C$, then there exists  $P_{\ell}$,  $P_{\ell-1}$ and  $P_{\ell-2}$ connecting $u,v$.
\end{prop}

\noindent{\bf Proof. }(1) Since $\zeta^{(\ell,k)}[B \cup C]=K_{\lfloor \frac{\ell-4}{2} \rfloor \, , \, \lfloor \frac{\ell-4}{2} \rfloor\,  ,  \, \left( \ell-4 - 2 \lfloor \frac{\ell-4}{2} \rfloor \right) , 1,1 }$ is a complete $r$-partite graph on $\ell-2$ vertices with $r \in \{3,4,5\}$, it is easy to check that any pairs of distinct vertices in $\zeta^{(\ell,k)}[B \cup C]$ are connected by  $P_{\ell-2}$, $P_{\ell-3}$ and  $P_{\ell-4}$.

	(2) By (1) there exists a $P_{\ell-2}$ and a $P_{\ell-3}$ connecting $b_4,b_5$ not use the vertices in $A$, so there  exists  $P_{\ell-1}$ and $P_{\ell-2}$ connecting $u,v$.
	
	(3) If $v \in A$, by (1) there exists $P_{\ell-2}$,  $P_{\ell-3}$ and $P_{\ell-4}$ connecting $b_4,b_5$ not use the vertices in $A$, so the result holds. If $v \in C$, by (1) there exists a $P_{\ell-2}$ and a $P_{\ell-3}$ connecting $b_4,v$ not use the vertices in $A$, so there is a $P_{\ell-1}$ and a $P_{\ell-2}$ connecting $u,v$. Moreover, consider the $P_{\ell-3}$ in $\zeta^{(\ell,k)}[B \cup C]$ connecting $b_4,v$, it is easy to see that such a $P_{\ell-3}$ can avoid $b_5$ by our construction. So by combining $ub_5w$ for some distinct $w \in A$ we can find a $P_{\ell}$ connecting $u,v$. \qed
\vspace{0.5em}

\noindent {\bf Remark.} Proposition \ref{P36} also holds for $\ell=5$ except for the $P_{\ell-4}=P_1$ connecting $u,v$ in (1). In the following paper we default that  Proposition \ref{P36} holds for all $\ell \geq 5$. ( We would not use the property that $u,v$ is connected by a $P_{\ell-4}$ in $\zeta^{(\ell,k)}[B \cup C]$ )

\begin{lemma}\label{L37}
For $k \geq \ell \geq 5$, $\zeta^{(\ell,k)}$ is $C_{\ell}$-saturated.
\end{lemma}

\noindent{\bf Proof. }We first prove that $\zeta^{(\ell,k)}$ is $C_{\ell}$-free. Using the same notations in Definition \ref{D35}, if $\zeta^{(\ell,k)}$ contains a $C_{\ell}$, then such a $C_{\ell}$ would contain at least two vertices in $A$. But all vertices in $A$ are only connected to $b_4,b_5$ which leads to a contradiction.
	
	Let $uv$ be an admissible nonedge in $\zeta^{(\ell,k)}$, then either $u \in A, v \in C$ or $u,v \in A$ by the construction. By Proposition \ref{P36} there exists a $P_{\ell}$ connecting $u,v$, so $\zeta^{(\ell,k)}+uv$ would create a $C_{\ell}$. Hence $\zeta^{(\ell,k)}$ is $C_{\ell}$-saturated. \qed
\vspace{0.5em}

	Similar to Definition \ref{D34}, we can construct a balanced $C_{\ell}$-saturated $k$-partite graph based on Proposition \ref{P31}.
\vspace{0.2em}

\begin{definition}\label{D38}
For $k \geq \ell \geq 5$ and $n \geq \lceil \frac{k+\ell-2}{2} \rceil$, let $Z^{(\ell,k,n)}$ be a graph   satisfies the following properties:  \\
\indent \text{(i)} $Z^{(\ell,k,n)}$ is a $k$-partite graph such that each part of $Z^{(\ell,k,n)}$ has exactly $n$ vertices. \\
\indent \text{(ii)} $Z^{(\ell,k,n)}=D_1 \cup D_2$ with $D_1 \cap D_2 = \emptyset$, $Z^{(\ell,k,n)}[D_1]=\zeta^{(\ell,k)}$  and $Z^{(\ell,k,n)}[D_2]$ is empty. \\
\indent \text{(iii)} Partition $D_1$ into $A,B,C$ with same notations to the definition of $\zeta^{(\ell,k)}$. Clearly $A$ is a good pair of $Z^{(\ell,k,n)}$.\\
\indent \text{(iv)} For any vertex $v \in D_2$, $d_{Z^{(\ell,k,n)}}(v)=1$, $N(v) \subseteq  A $. For different $u,v \in D_2$, if $u,v$ lie in the same part of $Z^{(\ell,k,n)}$, then $N(u)=N(v)$; otherwise  $N(u) \cap N(v) = \emptyset$.
\end{definition}
\vspace{0.0em}

	It is easy to check that $Z^{(\ell,k,n)}$ satisfies Proposition \ref{P31} by Proposition \ref{P36} and Lemma \ref{L37}. Lemma \ref{L33} guarantees the existence of $Z^{(\ell,k,n)}$. Let $\Xi^{(\ell,k,n)}$ be the graph family of all $Z^{(\ell,k,n)}$.
\vspace{0.8em}

\noindent{\bf Proof of Theorem \ref{T14}. }For $k\geq \ell \geq 5$ and $n \geq \lceil \frac{k+\ell-2}{2} \rceil$, let $G \in \Xi^{(\ell,k,n)}$, then $G$ is a $C_{\ell}$-saturated $k$-partite graph which implies 
\begin{equation*}
\begin{aligned}
sat(K_k^n, C_{\ell}) \leq e(G) & = kn+k+ \ell-5 + \lfloor \frac{\ell-4}{2} \rfloor ^2  + 2 \left( \ell-4-2  \lfloor \frac{\ell-4}{2} \rfloor \right) \lfloor \frac{\ell-4}{2} 
\rfloor \, .
\end{aligned}
\end{equation*}
   \qed

\subsection{Construction III}

In order to prove Theorem \ref{T110}, we construct a graph based on Proposition \ref{P31}. For $k \geq 3$ and $n,\ell$ sufficiently large, let $G$ be a $k$-partite graph with $n$ vertices in each part. Note that if $G$ satisfies the first two conditions of Proposition \ref{P31} but $G[A]$ is only $C_{\ell}$-free, then we can add some edges into $G[A]$ until $G[A]$ is $C_{\ell}$-saturated. Let $k \geq 4$ be a fixed integer.

	Let $C$ be a cycle of length $\ell-1$. We pick $x_1,\ldots,x_k$ along the cycle such that there are $s_i$ vertices between $x_i$ and $x_{i+1}$ (in the sense of module $k$), where $s_i \in \{\, \lfloor  \frac{\ell-1}{k} \rfloor -1  , \lceil \frac{\ell-1}{k} \rceil - 1  \} $ and $\sum \limits_{i=1}^k s_i = \ell-1-k$ (Figure 7). For $\ell \geq 60k+12$, let $m = \lceil \frac{\ell}{5} \rceil $ and $m_1,\ldots,m_k \in \{\,  \lfloor \frac{m}{k} \rfloor   , \lceil \frac{m}{k} \rceil \, \}$ such that $\sum \limits_{i=1}^k m_i = m$.

\begin{figure}[H]
\centering
\begin{tikzpicture}[scale=.45]

\draw (4,0) arc(0:360:4cm) ;

\filldraw (-2,3.45) circle (4pt);
\filldraw (2,3.45) circle (4pt);
\filldraw (3.9,0.95) circle (4pt);
\filldraw (0,-4) circle (4pt);
\filldraw (-3.9,0.95) circle (4pt);

\draw (-2,4.1) node[align=center]{$x_1$};
\draw (2,4.1) node[align=center]{$x_2$};
\draw (4.5,1.3) node[align=center]{$x_3$};
\draw (-4.5,1.3) node[align=center]{$x_k$};
\draw (0,-4.9) node[align=center]{$x_i$};

\draw[blue,thick](0,4.6) node[align=center]{$s_1$};
\draw[blue,thick](3.7,2.8) node[align=center]{$s_2$};
\draw[blue,thick](-3.7,2.8) node[align=center]{$s_k$};

\draw[dotted,thick] (3,0) arc(360:315:3cm) ;
\draw[dotted,thick] (-3,-0) arc(180:225:3cm) ;
\draw[dotted,thick] (3.5,-3.5) arc(315:355:5cm) ;
\draw[dotted,thick] (-3.5,-3.5) arc(225:185:5cm) ;

\draw[blue,thick] (-1.8,-4.3) node[align=center]{$s_{i}$};
\draw [blue,thick](2.1,-4.3) node[align=center]{$s_{i-1}$};

\end{tikzpicture}\\

\caption{\centering the vertices $\{x_1,\ldots,x_k\}$ on $C$.}

\end{figure}
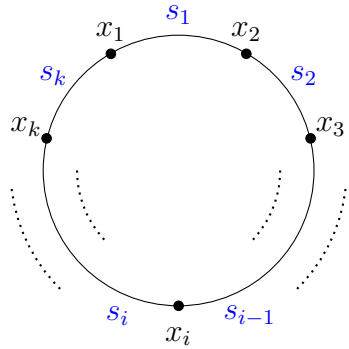

	A unit is obtained by adding 5 vertices and 9 edges on four consecutive vertices on $C$ (Figure 8). Let $P^i_{s_i+2}$ be the path connecting $x_i$ and $x_{i+1}$ of length $s_i+2$ on $C$, $1\le i\le k+1$. Assume that $P^i_{s_i+2}=x_i v_{i}^{1} v_{i}^{2} \ldots v_{i}^{s_i-1} v_{i}^{s_i} x_{i+1}$. We then add $m_i$ consecutive units into $P^i_{s_i+2}$ started from $v_i^4$ to $v_i^{4m_i+3}$ (Figure 9). Then there are still at least 3 vertices from $v_i^{4m_i+4}$ to $v_i^{s_i}$ since $s_i - (4m_i + 3 ) \geq \lfloor  \frac{\ell-1}{k} \rfloor -1 - (4  \lceil \lceil \frac{\ell}{5} \rceil / k \rceil +3 ) \geq 3$ while $\ell \geq 60k+12$. Let $H$ be the graph obtained by adding such units into $C$, $A=V(C)$ and $B=V(H) \setminus A$. Obviously, $H$ is a $k$-partite graph. Pick the Hamilton cycle $C^*$ in $H$ as Figure 10 shows.
\begin{figure}[H]
\centering

\begin{tikzpicture}[scale=0.8]

\draw[thick] (-2,1) -- (-1,1) -- (0,1) -- (1,1) -- (2,1) -- (3,1);

\draw[blue,thick] (-1,1) -- (-1,2) -- (0,2) -- (0,1) (1,1) -- (1,2) -- (2,2) -- (2,1);
\draw[blue,thick] (1,2) -- (1.5,2.7) -- (2,2)  (1.5,2.7) -- (2,1);

\filldraw (-1,1) circle (2pt);
\filldraw (0,1) circle (2pt);
\filldraw (1,1) circle (2pt);
\filldraw (2,1) circle (2pt);

\filldraw[red] (-1,2) circle (2pt);
\filldraw[red] (0,2) circle (2pt);
\filldraw[red] (1,2) circle (2pt);
\filldraw[red] (2,2) circle (2pt);
\filldraw[red] (1.5,2.7) circle (2pt);

\draw[thick](-2.4,1) node[align=center]{\large $C$};

\end{tikzpicture}\\

\caption{\centering A unit. The black lines and points represent edges and vertices on $C$. The blue lines and red points represent the addition of new edges and vertices. The two red vertices on the left are called unit-L, and the three red vertices on the right are called unit-R.}

\end{figure}
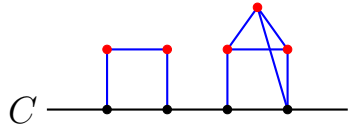

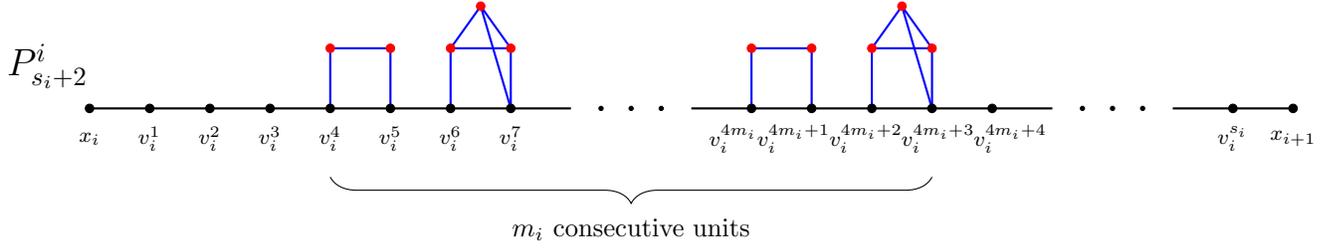
\begin{figure}[H]
\centering

\begin{tikzpicture}[scale=0.8]

\draw[thick] (-1,1) -- (7,1) ;

\filldraw (7.5,1) circle (1pt) (8,1) circle (1pt) (8.5,1) circle (1pt);

\filldraw (15.5,1) circle (1pt) (16,1) circle (1pt) (16.5,1) circle (1pt);

\draw[thick] (9,1) -- (10,1) -- (11,1) -- (12,1) -- (13,1) -- (14,1) -- (15,1)  (17,1) -- (19,1) ;

\draw[blue,thick] (3,1) -- (3,2) -- (4,2) -- (4,1) (5,1) -- (5,2) -- (6,2) -- (6,1);
\draw[blue,thick] (5,2) -- (5.5,2.7) -- (6,2)  (5.5,2.7) -- (6,1);

\draw[blue,thick] (10,1) -- (10,2) -- (11,2) -- (11,1) (12,1) -- (12,2) -- (13,2) -- (13,1);
\draw[blue,thick] (12,2) -- (12.5,2.7) -- (13,2)  (12.5,2.7) -- (13,1);

\filldraw (-1,1) circle (2pt);
\filldraw (0,1) circle (2pt);
\filldraw (1,1) circle (2pt);
\filldraw (2,1) circle (2pt);

\filldraw (3,1) circle (2pt);
\filldraw (4,1) circle (2pt);
\filldraw (5,1) circle (2pt);
\filldraw (6,1) circle (2pt);

\filldraw[red] (3,2) circle (2pt);
\filldraw[red] (4,2) circle (2pt);
\filldraw[red] (5,2) circle (2pt);
\filldraw[red] (6,2) circle (2pt);
\filldraw[red] (5.5,2.7) circle (2pt);

\filldraw (10,1) circle (2pt);
\filldraw (11,1) circle (2pt);
\filldraw (12,1) circle (2pt);
\filldraw (13,1) circle (2pt);

\filldraw[red] (10,2) circle (2pt);
\filldraw[red] (11,2) circle (2pt);
\filldraw[red] (12,2) circle (2pt);
\filldraw[red] (13,2) circle (2pt);
\filldraw[red] (12.5,2.7) circle (2pt);

\filldraw (14,1) circle (2pt);
\filldraw (18,1) circle (2pt);
\filldraw (19,1) circle (2pt);

\draw[thick](-1.7,1.7) node[align=center]{\large $P^i_{s_i+2}$};

\draw (-1,0.5) node[align=center]{\scriptsize $x_i$};
\draw (0,0.5) node[align=center]{\scriptsize $v_{i}^{1}$};
\draw (1,0.5) node[align=center]{\scriptsize $v_{i}^{2}$};
\draw (2,0.5) node[align=center]{\scriptsize $v_{i}^{3}$};
\draw (3,0.5) node[align=center]{\scriptsize $v_{i}^{4}$};
\draw (4,0.5) node[align=center]{\scriptsize $v_{i}^{5}$};
\draw (5,0.5) node[align=center]{\scriptsize $v_{i}^{6}$};
\draw (6,0.5) node[align=center]{\scriptsize $v_{i}^{7}$};

\draw [decorate,decoration={brace,amplitude=10pt,mirror},yshift=-4pt](3,0) -- (13,0) node[black,midway,yshift=-0.7cm] {\footnotesize $m_i$ consecutive units};

\draw (9.7,0.5) node[align=center]{\scriptsize $v_{i}^{4m_i}$};
\draw (10.7,0.5) node[align=center]{\scriptsize $v_{i}^{4m_i+1}$};
\draw (11.9,0.5) node[align=center]{\scriptsize $v_{i}^{4m_i+2}$};
\draw (13.1,0.5) node[align=center]{\scriptsize $v_{i}^{4m_i+3}$};
\draw (14.3,0.5) node[align=center]{\scriptsize $v_{i}^{4m_i+4}$};

\draw (18,0.5) node[align=center]{\scriptsize $v_{i}^{s_i}$};
\draw (19,0.5) node[align=center]{\scriptsize $x_{i+1}$};

\end{tikzpicture}\\

\caption{\centering The addition of $m_i$ units into $P^i_{s_i+2}$.}

\end{figure}

\begin{figure}[H]
\centering

\begin{tikzpicture}[scale=0.8]

\draw[thick] (-1,1) -- (3,1) (4,1) -- (5,1) (6,1) -- (7,1);

\filldraw (7.5,1) circle (1pt) (8,1) circle (1pt) (8.5,1) circle (1pt);

\filldraw (15.5,1) circle (1pt) (16,1) circle (1pt) (16.5,1) circle (1pt);

\draw[thick] (9,1) -- (10,1) (11,1) -- (12,1)  (13,1) -- (14,1) -- (15,1)  (17,1) -- (19,1) ;

\draw[thick] (3,1) -- (3,2) -- (4,2) -- (4,1) (5,1) -- (5,2)  (6,2) -- (6,1);
\draw[thick] (5,2) -- (5.5,2.7) -- (6,2)  ;

\draw[thick] (10,1) -- (10,2) -- (11,2) -- (11,1) (12,1) -- (12,2)  (13,2) -- (13,1);
\draw[thick] (12,2) -- (12.5,2.7) -- (13,2)  ;

\filldraw (-1,1) circle (2pt);
\filldraw (0,1) circle (2pt);
\filldraw (1,1) circle (2pt);
\filldraw (2,1) circle (2pt);

\filldraw (3,1) circle (2pt);
\filldraw (4,1) circle (2pt);
\filldraw (5,1) circle (2pt);
\filldraw (6,1) circle (2pt);

\filldraw (3,2) circle (2pt);
\filldraw (4,2) circle (2pt);
\filldraw (5,2) circle (2pt);
\filldraw (6,2) circle (2pt);
\filldraw (5.5,2.7) circle (2pt);

\filldraw (10,1) circle (2pt);
\filldraw (11,1) circle (2pt);
\filldraw (12,1) circle (2pt);
\filldraw (13,1) circle (2pt);

\filldraw (10,2) circle (2pt);
\filldraw (11,2) circle (2pt);
\filldraw (12,2) circle (2pt);
\filldraw (13,2) circle (2pt);
\filldraw (12.5,2.7) circle (2pt);

\filldraw (14,1) circle (2pt);
\filldraw (18,1) circle (2pt);
\filldraw (19,1) circle (2pt);

\draw[thick](-1.7,1.7) node[align=center]{\large $P^i_{s_i+2}$};

\draw (-1,0.5) node[align=center]{\scriptsize $x_i$};
\draw (19,0.5) node[align=center]{\scriptsize $x_{i+1}$};

\end{tikzpicture}\\

\caption{\centering The edges of Hamilton cycle $C^*$ in $P^i_{s_i+2}$.}

\end{figure}
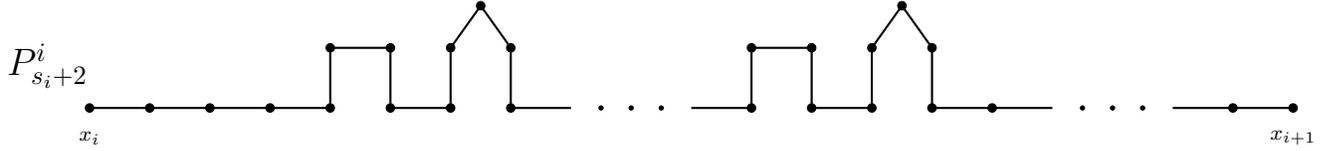

Let $H=V'_1\cup\dots\cup V'_k$. In order to construct a balanced $C_{\ell}$-saturated $k$-partite graph, we mark some special vertices. See Figure 11, for any $i\in\{1,\ldots,k\}$,  let $v_i^{2}, v_{i-1}^{s_i-1} \in V'_i$ (including all vertices with the same color as them) and $x_i \not \in V'_i$. This can be done since $k \geq 4$. Now we add some edges into $H$ (see Figure 11). A  detailed explanation in adding the edges can be find in Figure 12. We  denote the obtained graph by $H'$. 

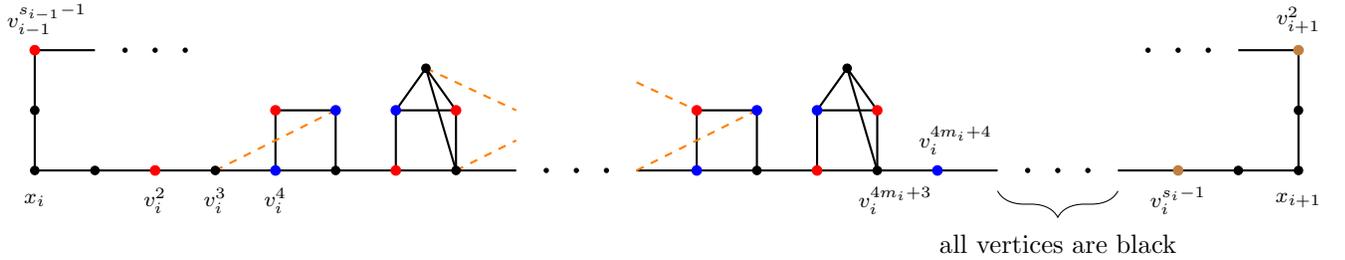
\begin{figure}[H]
\centering

\begin{tikzpicture}[scale=0.8]

\draw[thick] (0,3) -- (-1,3) -- (-1,1) -- (7,1) ;

\filldraw (7.5,1) circle (1pt) (8,1) circle (1pt) (8.5,1) circle (1pt);

\filldraw (15.5,1) circle (1pt) (16,1) circle (1pt) (16.5,1) circle (1pt);

\filldraw (0.5,3) circle (1pt) (1,3) circle (1pt) (1.5,3) circle (1pt);

\filldraw (18.5,3) circle (1pt) (18,3) circle (1pt) (17.5,3) circle (1pt);

\draw[thick] (9,1) -- (10,1) -- (11,1) -- (12,1) -- (13,1) -- (14,1) -- (15,1)  (17,1) -- (20,1) -- (20,3) -- (19,3);

\draw[thick] (3,1) -- (3,2) -- (4,2) -- (4,1) (5,1) -- (5,2) -- (6,2) -- (6,1);
\draw[thick] (5,2) -- (5.5,2.7) -- (6,2)  (5.5,2.7) -- (6,1);

\draw[thick] (10,1) -- (10,2) -- (11,2) -- (11,1) (12,1) -- (12,2) -- (13,2) -- (13,1);
\draw[thick] (12,2) -- (12.5,2.7) -- (13,2)  (12.5,2.7) -- (13,1);

\draw[orange,dashed,thick] (2,1) -- (4,2)   (6,1) -- (7,1.5)  (5.5,2.7) -- (7,2);

\draw[orange,dashed,thick] (9,1) -- (11,2)  (9,2.466) -- (10,2);

\filldraw[red,thick] (-1,3) circle (2pt);
\filldraw (-1,2) circle (2pt);
\filldraw (-1,1) circle (2pt);
\filldraw (0,1) circle (2pt);
\filldraw[red,thick] (1,1) circle (2pt);
\filldraw (2,1) circle (2pt);

\filldraw[blue,thick] (3,1) circle (2pt);
\filldraw (4,1) circle (2pt);
\filldraw[red,thick] (5,1) circle (2pt);
\filldraw (6,1) circle (2pt);

\filldraw[red,thick] (3,2) circle (2pt);
\filldraw[blue,thick] (4,2) circle (2pt);
\filldraw[blue,thick] (5,2) circle (2pt);
\filldraw[red,thick] (6,2) circle (2pt);
\filldraw (5.5,2.7) circle (2pt);

\filldraw[blue,thick] (10,1) circle (2pt);
\filldraw (11,1) circle (2pt);
\filldraw[red,thick] (12,1) circle (2pt);
\filldraw (13,1) circle (2pt);

\filldraw[red,thick] (10,2) circle (2pt);
\filldraw[blue,thick] (11,2) circle (2pt);
\filldraw[blue,thick] (12,2) circle (2pt);
\filldraw[red,thick] (13,2) circle (2pt);
\filldraw (12.5,2.7) circle (2pt);

\filldraw[blue,thick] (14,1) circle (2pt);
\filldraw[brown,thick] (18,1) circle (2pt);
\filldraw (19,1) circle (2pt);
\filldraw (20,1) circle (2pt);
\filldraw (20,2) circle (2pt);
\filldraw[brown,thick] (20,3) circle (2pt);

\draw [decorate,decoration={brace,amplitude=10pt,mirror},yshift=-4pt](15,0.8) -- (17,0.8) node[black,midway,yshift=-0.7cm] {\footnotesize all vertices are black};

\draw (-0.8,3.5) node[align=center]{\scriptsize $v_{i-1}^{s_{i-1}-1}$};
\draw (-1,0.5) node[align=center]{\scriptsize $x_i$};
\draw (20,0.5) node[align=center]{\scriptsize $x_{i+1}$};
\draw (20,3.5) node[align=center]{\scriptsize $v_{i+1}^{2}$};
\draw (1,0.5) node[align=center]{\scriptsize $v_{i}^{2}$};
\draw (2,0.5) node[align=center]{\scriptsize $v_{i}^{3}$};
\draw (3,0.5) node[align=center]{\scriptsize $v_{i}^{4}$};
\draw (13.3,0.5) node[align=center]{\scriptsize $v_{i}^{4m_i+3}$};
\draw (14.3,1.5) node[align=center]{\scriptsize $v_{i}^{4m_i+4}$};

\draw (18,0.5) node[align=center]{\scriptsize $v_{i}^{s_i-1}$};

\end{tikzpicture}\\

\caption{\centering The addition of edges and the partition of vertices in $P^i_{s_i+2}$. The same color (red, blue, brown) of vertices are in the same part of $H'$ except black color. The black vertices can lie in any part of $H'$ different from the neighbours. The orange dashed lines are the additional edges and the detail would be shown in the next figure.}

\end{figure}

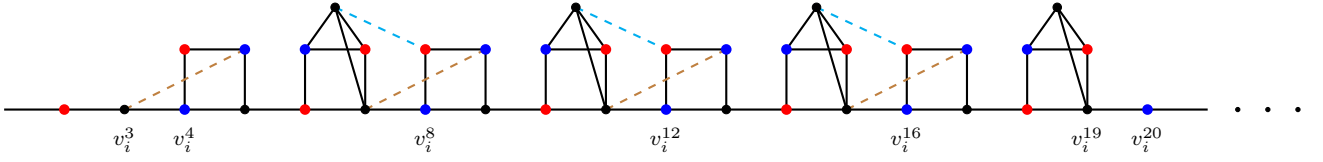
\begin{figure}[H]
\centering

\begin{tikzpicture}[scale=0.8]

\draw[thick] (0,1) -- (7,1) ;

\filldraw (20.5,1) circle (1pt) (21,1) circle (1pt) (21.5,1) circle (1pt);

\draw[thick] (7,1) -- (10,1) -- (11,1) -- (12,1) -- (13,1) -- (14,1) -- (20,1)  ;

\draw[thick] (3,1) -- (3,2) -- (4,2) -- (4,1) (5,1) -- (5,2) -- (6,2) -- (6,1);
\draw[thick] (5,2) -- (5.5,2.7) -- (6,2)  (5.5,2.7) -- (6,1);

\draw[thick] (7,1) -- (7,2) -- (8,2) -- (8,1) (9,1) -- (9,2) -- (10,2) -- (10,1);
\draw[thick] (9,2) -- (9.5,2.7) -- (10,2)  (9.5,2.7) -- (10,1);

\draw[thick] (11,1) -- (11,2) -- (12,2) -- (12,1) (13,1) -- (13,2) -- (14,2) -- (14,1);
\draw[thick] (13,2) -- (13.5,2.7) -- (14,2)  (13.5,2.7) -- (14,1);

\draw[thick] (15,1) -- (15,2) -- (16,2) -- (16,1) (17,1) -- (17,2) -- (18,2) -- (18,1);
\draw[thick] (17,2) -- (17.5,2.7) -- (18,2)  (17.5,2.7) -- (18,1);

\draw[brown,dashed,thick] (2,1) -- (4,2)   (6,1) -- (8,2)  ;
\draw[cyan,dashed,thick] (5.5,2.7) -- (7,2) (9.5,2.7) -- (11,2) (13.5,2.7) -- (15,2) ;

\draw[brown,dashed,thick] (10,1) -- (12,2)   (14,1) -- (16,2)    ;

\filldraw[red,thick] (1,1) circle (2pt);
\filldraw (2,1) circle (2pt);

\filldraw[blue,thick] (3,1) circle (2pt);
\filldraw (4,1) circle (2pt);
\filldraw[red,thick] (5,1) circle (2pt);
\filldraw (6,1) circle (2pt);

\filldraw[red,thick] (3,2) circle (2pt);
\filldraw[blue,thick] (4,2) circle (2pt);
\filldraw[blue,thick] (5,2) circle (2pt);
\filldraw[red,thick] (6,2) circle (2pt);
\filldraw (5.5,2.7) circle (2pt);

\filldraw[blue,thick] (7,1) circle (2pt);
\filldraw (8,1) circle (2pt);
\filldraw[red,thick] (9,1) circle (2pt);
\filldraw (10,1) circle (2pt);

\filldraw[red,thick] (7,2) circle (2pt);
\filldraw[blue,thick] (8,2) circle (2pt);
\filldraw[blue,thick] (9,2) circle (2pt);
\filldraw[red,thick] (10,2) circle (2pt);
\filldraw (9.5,2.7) circle (2pt);

\filldraw[blue,thick] (11,1) circle (2pt);
\filldraw (12,1) circle (2pt);
\filldraw[red,thick] (13,1) circle (2pt);
\filldraw (14,1) circle (2pt);

\filldraw[red,thick] (11,2) circle (2pt);
\filldraw[blue,thick] (12,2) circle (2pt);
\filldraw[blue,thick] (13,2) circle (2pt);
\filldraw[red,thick] (14,2) circle (2pt);
\filldraw (13.5,2.7) circle (2pt);

\filldraw[blue,thick] (15,1) circle (2pt);
\filldraw (16,1) circle (2pt);
\filldraw[red,thick] (17,1) circle (2pt);
\filldraw (18,1) circle (2pt);

\filldraw[red,thick] (15,2) circle (2pt);
\filldraw[blue,thick] (16,2) circle (2pt);
\filldraw[blue,thick] (17,2) circle (2pt);
\filldraw[red,thick] (18,2) circle (2pt);
\filldraw (17.5,2.7) circle (2pt);

\filldraw[blue,thick] (19,1) circle (2pt);
\draw (2,0.5) node[align=center]{\scriptsize $v_{i}^{3}$};
\draw (3,0.5) node[align=center]{\scriptsize $v_{i}^{4}$};
\draw (7,0.5) node[align=center]{\scriptsize $v_{i}^{8}$};
\draw (11,0.5) node[align=center]{\scriptsize $v_{i}^{12}$};
\draw (15,0.5) node[align=center]{\scriptsize $v_{i}^{16}$};
\draw (18,0.5) node[align=center]{\scriptsize $v_{i}^{19}$};
\draw (19,0.5) node[align=center]{\scriptsize $v_{i}^{20}$};

\end{tikzpicture}\\

\caption{\centering The addition of edges between units (detailed figure). This is an example on exactly $m_i=4$ units. There are two kinds of edges, cyan and brown. Cyan edges $uv$ satisfy $u,v \in B$. Brown edges $uv$ satisfy $u \in A$ and $v \in B$. For $m_i$ units there are $m_i-1$ cyan edges and $m_i$ brown edges.}

\end{figure}

The construction of $H'$ gives us some properties.
\vspace{0.1em} 	

\begin{lemma}\label{L39}  $H'$ is a $k$-partite $C_{\ell}$-free graph.
\end{lemma}
\noindent{\bf Proof. }By the construction of $H'$ and $k\ge 4$, it is not difficult to check that $H'$ is a $k$-partite graph.
Let $uv\in E(H')$. We will show that there is no path of length $\ell$ connecting $u$ and $v$, then there is no  $C_{\ell}$ in $H'$. Recall $\ell \geq 60k+12$. In the following discussion, we only mention the paths of lengths at least $\ell-1$ and ignore short ones.

If $u,v \in A$, by the construction of $H'$, the lengths of paths connecting $u$ and $v$ in $H'-uv$ are $s$ with  $s=\ell-1$ or $s\ge \ell+1$.  If $u \in B $, $v \in A$ and $uv$ is not a brown edge or $u,v \in B$ in Figure 12, then the length of the shortest path in $H'-uv$  connecting $u$ and $v$ is at least $\ell+1$ by the construction of $H'$.

If $uv$ is a brown edge in Figure 12, by the construction of $H'$, the lengths of paths connecting $u$ and $v$ in $H'-uv$ are $s$ with  $s=\ell-1$ or $s\ge \ell+1$.
Thus we are done.\qed
\vspace{0.3em}

For any $u,v \in A$, let $d_C^-(u,v)$ be the distance between $u,v$ on $C$ and $d_C^+(u,v)= \ell-1- d_C^-(u,v) $. Clearly, $0 \leq d_C^-(u,v) \leq \frac{\ell-1}{2} \leq d_C^+(u,v) \leq \ell-1$ and there exists a path of length $d_C^+(u,v)+1$ on $C$ connecting $u$ and $v$. If $u=v$ we define $d_C^-(u,v)=0$ and $d_C^+(u,v) =\ell-1$. For any $u,v \in V(H')$, let $d^{\max}(u,v)$ be the number of edges in the longest path connecting $u$ and $v$ in $H'$. Since $m = \lceil \frac{\ell}{5} \rceil $, we have $|H'| = \ell + 5m \geq 2 \ell $. Thus $d^{\max}(u,v) \geq \ell$ because there exists a $P_r$ connecting $u$ and $v$ in $C^*$ where $r \geq \ell+1$. For any $u \in B$, let $u^A$ be the only neighbour of $u$ in $A$. For a pair of neighbours $u,v \in A$, if the neighbours of $u,v$ in $B$ lie in a same unit-L or unit-R, we say $u,v$ are \textit{basic neighbours}.
	
	For any different $u,v \in V(H')$, let 
\begin{equation*}
d^\min (u,v) =
\left\{ \;
\begin{aligned}
& d^+_C (u,v) \, , & \; & u,v \in A   \, ; \\
& d^+_C (u^A,v)+1 \, , & \; & u \in B, v \in A  \, ; \\
& d^+_C (u^A,v^A)+2 \, , & \; & u,v \in B \, .
\end{aligned}
\right.
\end{equation*}

 By the construction of $H'$, $d^{\max}(u,v)\ge d^\min (u,v)+2$ for any different $u,v \in V(H')$.
\vspace{0.1em}

\begin{lemma}\label{L310}  For any different $u,v \in V(H')$ and $q \in [ \, d^\min (u,v)+2,d^\max (u,v) \, ]\cup\{d^\min (u,v)\}$, there exists a path of length $q+1$ connecting $u$ and $v$ in $H'$.
\end{lemma}
\noindent{\bf Proof. }Let $u,v\in V(H')$ with $u\not=v$. We first show that there is a path, denoted by $P^\min(u,v)$,  of length $d^\min (u,v)+1$ connecting $u$ and $v$ in $H'$. If $u,v \in A$, then  $P^\min(u,v)$ is the longest path on $C$ connecting $u$ and $v$.
If $u \in B, v \in A$ and $u^A \neq v$, then $P^\min(u,v) = uP^\min(u^A,v)$. If $u \in B, v \in A$ and $u^A = v$, then $d^+_C (u^A,v) = \ell-1$. Let $w$ be the basic neighbour of $v$ and $z$ be the vertex such that $z^A=w$ and $uz \in E(H')$. Then $P^\min(u,v) = uzP^\min(w,v)$ who has length $\ell+1$.
	If $u,v \in B$ and $u^A \neq v^A$ then $P^\min(u,v) = uP^\min(u^A,v^A)v$. If $u,v \in B$ and $u^A = v^A$ then $u ,v$ must lie in some same unit-R and $d^+_C (u^A,v^A) = \ell-1$. Let $w$ be the basic neighbour of $u^A$ and $z$ be the rest vertex in the same unit-R. Then $P^\min(u,v) = uP^\min(u^A,w)zv$ who has length $\ell+2$. In both cases, the length of $P^\min(u,v)$ is $d^\min (u,v)+1$.

Now we  expand the path $P^\min(u,v)$ by adding some units into it, where the edges are chosen blue ones in Figure 8. Since a unit can be divided into two parts (\,2-vertices unit-L and 3-vertices unit-R\,), we can choose 2 vertices in unit-L and either 2 or 3 vertices in unit-R to expand the path. By B\'ezout's Lemma, for any integer $q \in [ \, d^\min (u,v)+2,d^\max (u,v) \, ]$,
we can find suitable combination of unit-L's and unit-R's to expand the path to length $q+1$. \qed
\vspace{0.1em}

\begin{cor}\label{C311}
For fixed $x_i$ and any different $u \in V(H')$, we have\\
\indent (1) There exists a $P_{\ell-2}$ connecting $x_i$ and $x_j$, where $j \neq i$; \\
\indent (2) there exists a $P_{\ell-1}$ connecting $x_i$ and $u$ except $d(u,x_i)=2$.
\end{cor}
\noindent{\bf Proof. }(1)\;It is easy to check that $ \frac{\ell-1}{2} \leq d^\min(x_i,x_j) \leq \ell-1- \lfloor \frac{\ell-1}{k} \rfloor$ and $d^\max(x_i,x_j) \geq \ell$. By Lemma \ref{L310} we immediately have the result. \\
\indent (2)\;If $d(x_i,u) \geq 3$ then $d^\min(x_i,u) \leq \ell-4$. By Lemma \ref{L310}, there exists a path with length $\ell-1$ connecting $x_i$ and $u$. If $d(x_i,u)=1$ then $P^\min(x_i,u)$ is the $P_{\ell-1}$ connecting them. \qed
\vspace{0.5em}

\begin{lemma}\label{L312} $H'$ is a $k$-partite $C_{\ell}$-saturated graph.
\end{lemma}
\noindent{\bf Proof. } By Lemma \ref{L39}, we just need to show that for any different $u,v \in V(H')$ and $uv$ being an admissible nonedge,  there is a $P_{\ell}$ connecting $u$ and $v$  and then there is a $C_{\ell}$ in $H'+uv$.
By Lemma \ref{L310}, it suffices to show that $d^\min (u,v)\le \ell-3$ or $d^\min (u,v)= \ell-1$ or there is path of length $ \ell$  connecting $u$ and $v$. We also use $P^\min(u,v)$ to denote the path  of length $d^\min (u,v)+1$ connecting $u$ and $v$ in $H'$.

If $u,v \in A$, then $d^-_C(u,v)\ge 2$ which implies $d^+_C(u,v)\le \ell-3$. Thus $d^\min (u,v)\le \ell-3$ and we are done.

Let $u,v\in B$. Since $uv$ is an admissible nonedge, by the construction of $H'$, $d^-_C (u^A,v^A)\ge 2$. If $d^-_C (u^A,v^A)\ge 4$ (resp. $d^-_C (u^A,v^A)= 2$), then $d^\min (u,v)\le \ell-3$ (resp. $d^\min (u,v)= \ell-1$) and we are done. So we just consider the case $d^-_C (u^A,v^A)=3$. Then we have $\{u,v\}\in\{\{z_1,w_1\},\{z_2,z_5\},\{w_1,z_7\}\}$ (see Figure 13). If $\{u,v\}=\{z_1,w_1\}$, say $u=z_1$ and $v=w_1$, then $uz_2P^\min(y_0,y_4)z_4v$ is a path connecting $u$ and $v$  of length $\ell$ and we are done. If $\{u,v\}=\{z_2,z_5\}$, say $u=z_2$ and $v=z_5$, then $uP^\min(y_0,y_8)z_8w_2z_7y_7y_6y_5v$ is a path connecting $u$ and $v$ of length $\ell$ and we are done. If $\{u,v\}=\{w_1,z_7\}$, say $u=w_1$ and $v=z_7$, then $uz_4P^\min(y_4,y_7)v$ is a path connecting $u$ and $v$ of length $\ell$ and we are done.

Let $u \in B$ and $v \in A$. If $d^-_C (u^A,v)\ge 3$ (resp. $d^-_C (u^A,v)= 1$), then $d^\min (u,v)\le \ell-3$ (resp. $d^\min (u,v)= \ell-1$) and we are done. So we just consider the case $d^-_C (u^A,v)=2$.  Then we have $u=z_2,v=y_4$ or $u\in\{w_1,z_4\},v\in\{y_2,y_6\}$  (see Figure 13). If $u=z_2,v=y_4$, then $uP^\min(y_0,y_8)z_8w_2z_7y_7y_6y_5v$ is a path connecting $u$ and $v$ of length $\ell$ and we are done. If $u\in\{w_1,z_4\},v\in\{y_2,y_6\}$, say $u=w_1$ and $v=y_2$, then $P^\min(v,y_4)z_4 u$ is a path connecting $u$ and $v$ of length $\ell$ and we are done.
 \qed
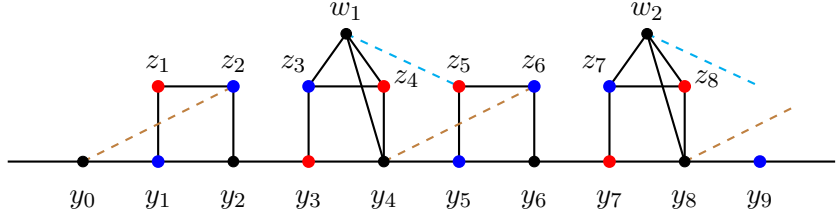
\begin{figure}[H]
\centering

\begin{tikzpicture}[scale=1]

\draw[thick] (1,1) -- (7,1) ;

\draw[thick] (7,1) -- (10,1) -- (11,1) -- (12,1)  ;

\draw[thick] (3,1) -- (3,2) -- (4,2) -- (4,1) (5,1) -- (5,2) -- (6,2) -- (6,1);
\draw[thick] (5,2) -- (5.5,2.7) -- (6,2)  (5.5,2.7) -- (6,1);

\draw[thick] (7,1) -- (7,2) -- (8,2) -- (8,1) (9,1) -- (9,2) -- (10,2) -- (10,1);
\draw[thick] (9,2) -- (9.5,2.7) -- (10,2)  (9.5,2.7) -- (10,1);


\draw[brown,dashed,thick] (2,1) -- (4,2)   (6,1) -- (8,2)  ;
\draw[cyan,dashed,thick] (5.5,2.7) -- (7,2) (9.5,2.7) -- (11,2) ;

\draw[brown,dashed,thick] (10,1) -- (11.5,1.75)   ;

\filldraw (2,1) circle (2pt);

\filldraw[blue,thick] (3,1) circle (2pt);
\filldraw (4,1) circle (2pt);
\filldraw[red,thick] (5,1) circle (2pt);
\filldraw (6,1) circle (2pt);

\filldraw[red,thick] (3,2) circle (2pt);
\filldraw[blue,thick] (4,2) circle (2pt);
\filldraw[blue,thick] (5,2) circle (2pt);
\filldraw[red,thick] (6,2) circle (2pt);
\filldraw (5.5,2.7) circle (2pt);

\filldraw[blue,thick] (7,1) circle (2pt);
\filldraw (8,1) circle (2pt);
\filldraw[red,thick] (9,1) circle (2pt);
\filldraw (10,1) circle (2pt);

\filldraw[red,thick] (7,2) circle (2pt);
\filldraw[blue,thick] (8,2) circle (2pt);
\filldraw[blue,thick] (9,2) circle (2pt);
\filldraw[red,thick] (10,2) circle (2pt);
\filldraw (9.5,2.7) circle (2pt);

\filldraw[blue,thick] (11,1) circle (2pt);

\draw (2,0.5) node[align=center]{\small $y_0$};
\draw (3,0.5) node[align=center]{\small $y_1$};
\draw (4,0.5) node[align=center]{\small $y_2$};
\draw (5,0.5) node[align=center]{\small $y_3$};
\draw (6,0.5) node[align=center]{\small $y_4$};
\draw (7,0.5) node[align=center]{\small $y_5$};
\draw (8,0.5) node[align=center]{\small $y_6$};
\draw (9,0.5) node[align=center]{\small $y_7$};
\draw (10,0.5) node[align=center]{\small $y_8$};
\draw (11,0.5) node[align=center]{\small $y_9$};

\draw (3,2.3) node[align=center]{\small $z_1$};
\draw (4,2.3) node[align=center]{\small $z_2$};
\draw (4.8,2.3) node[align=center]{\small $z_3$};
\draw (6.3,2.1) node[align=center]{\small $z_4$};
\draw (7,2.3) node[align=center]{\small $z_5$};
\draw (8,2.3) node[align=center]{\small $z_6$};
\draw (8.8,2.3) node[align=center]{\small $z_7$};
\draw (10.3,2.1) node[align=center]{\small $z_8$};

\draw (5.5,3) node[align=center]{\small $w_1$};
\draw (9.5,3) node[align=center]{\small $w_2$};

\end{tikzpicture}\\

\caption{Auxiliary figure for Lemma \ref{L312}.}

\end{figure}

\vspace{0.1em}

 	Now we give our construction.

\noindent \textbf{Construction of $G$. }For $k \geq 3$, $\ell \geq 60k+12$  and $n$ sufficiently large, let $G=V_1 \cup \cdots \cup V_k$ be the balanced $k$-partite graph with $n$ vertices in each part such that  $H'$ is a subgraph of $G$ with $V(H^')=V_1^{'}\cup V_2^{'}\cup ...\cup V_{k}^{'}$ where $V_i' \subseteq V_i$. For any $u_i \in V_i  \setminus V'_i \neq \emptyset$, let $N(u_i)=\{x_i\}$.   
\vspace{0.3em}

Let $\mathcal{A}=V(H')$. By Lemmas \ref{L39}, \ref{L312}, Corollary \ref{C311} and Proposition \ref{P31}, $G$ is a $k$-partite $C_{\ell}$-saturated graph. Now we can prove Theorem \ref{T110}.

\noindent{\bf Proof of Theorem \ref{T110}. }By the Construction of $G$, we have
\begin{equation*}
\begin{aligned}
sat(K_k^n,C_{\ell}) & \leq e(G)  = | \, G \setminus H' \, | \, + \, e(H') \\
& = kn-(\ell-1+5m) + \left( \ell-1 + 9m + \sum_{i=1}^k (m_i-1) + \sum_{i=1}^k m_i \right) \\
& = kn +6m - k \\
& = k(n-1) + 6 \, \lceil \,  \frac{\ell}{5}  \, \rceil \, .
\end{aligned}
\end{equation*}
The case $k=3$ will be given in the following Remark. \qed
\vspace{0.5em}

\noindent \textbf{Remark. }Although the construction above only gave the case $k \geq 4$, we can similarly construct the case $k =3$ since $H$ is actually a tripartite graph. The difference between them is the number of extra additional edges. The proof of Lemma \ref{L312} guarantees that no more than $9m$ cyan or brown edges would be added into $H$ to obtain $H'$, and we still have $sat(K_3^n,C_{\ell}) \leq 3n + 13m = 3n + 13 \lceil \,  \frac{\ell}{5}  \, \rceil $. \qed

\section{$C_4$-saturation multipartite graph}	

\subsection{Tripartite}

For a bipartite graph $H=A \cup B$, we say $H$ is \textit{perfect} to $A$ if for any vertex  $u\in A$, $d_H(u)\ge 1$; otherwise  $H$ is \textit{non-perfect} to $A$.

Let $G \in Sat(K_{3} ^{ n}, C_4)$, we write $G=V_1 \cup V_2 \cup V_3$. It is easy to check that $W^{(4,3,n)}$ (see Figure 2) is a $C_4$-saturated tripartite graph, thus $e(G) \leq 3n$. The following lemmas show that  $Sat(K_{3} ^{ n}, C_4)=\{W^{(4,3,n)}\}$.
\vspace{0.0em}

\begin{lemma}\label{L41}  For any distinct $i,j \in [3]$, $G \, [ \, V_i \cup V_j \, ] $ has no isolated vertices $a,b$ such that $a \in V_i$ and $b \in V_j$.
\end{lemma}

\noindent{\bf Proof. }Suppose that there exists $i,j \in [3]$ such that $G \, [ \, V_i \cup V_j \, ] $ has isolated vertices $a \in V_i$ and $b \in V_j$. Then  $G+ab$ would form a $C_4$, say $abuva$ in $G$.  Since $a,b$ are isolated vertices in $G \, [ \, V_i \cup V_j \, ] $,  $u,v \in V_k$ where $k \in [3] \setminus \{i,j\}$, a contradiction with
 $uv \not \in E(G)$. \qed
\vspace{0.4em}

	We immediately have the following corollaries by Lemma \ref{L41} and $e(G)\le 3n$.

\begin{cor}\label{C42}  For any distinct $i,j \in [3]$, $G \, [ \, V_i \cup V_j \, ] $ is perfect to either $V_i$ or $V_j$.
\end{cor}

\begin{cor}\label{C43}  For any distinct $i,j \in [3]$, $e \left( G \, [ \, V_i \cup V_j \, ] \right) = n $. %

\end{cor}

	The next two lemmas show the core properties of $G$.

\begin{lemma}\label{L44}  If $G \, [ \, V_i \cup V_j \, ] $ is a perfect matching for any distinct $i,j \in [3]$, then $G$ is not $C_4$-saturated.
\end{lemma}
\noindent{\bf Proof. }We construct a digraph $G'$ by giving all edges of $G$ a direction from $V_1 \rightarrow V_2 \rightarrow V_3 \rightarrow V_1$. Then for any $u \in V(G')$, $d^+(u)=d^-(u)=1$. Since $G$ is connected, $G'$ must be a directed Hamilton cycle which can be assumed as
\begin{equation*}
  x_1 \rightarrow y_1 \rightarrow z_1 \rightarrow x_2 \rightarrow y_2 \rightarrow z_2 \rightarrow \ldots \rightarrow x_n \rightarrow y_n \rightarrow z_n \rightarrow x_1
\end{equation*}
where $x_i \in V_1 , y_i \in V_2 , z_i \in V_3$. Then  $G+x_1y_2$ can only create a $C_5$ and a $C_{3n-3}$ in $G$. Thus $G$ is not $C_4$-saturated. \qed

\begin{lemma}\label{L45}  Let $i,j,k$ be distinct elements of $\{1,2,3\}$. If $G \, [ \, V_i \cup V_j \, ] $ is perfect to $V_i$ and non-perfect to $V_j$, then there exists a unique vertex $v\in V_k$ such that $N_j(v)=V_j$.
\end{lemma}

\noindent{\bf Proof. }Since $G \, [ \, V_i \cup V_j \, ] $ is perfect to $V_i$ and non-perfect to $V_j$, $N_j(u)\not=\emptyset$ for any $u \in V_i$ and $ V_j \setminus N_j(V_i)\not=\emptyset$. By  Corollary \ref{C43},  $|N_j(u)|=1$ for any $u \in V_i$.

If there exists $v_j\in V_j$ such that $|N_k(v_j)|=0$,  consider $G+v_jv_k$ for any $v_k\in V_k$. Suppose the $C_4$ in $G+v_jv_k$ is $v_jw_1w_2v_kv_j$, then $w_1\in V_i$ and $w_2\in V_j$. Hence $N_j(w_1)=2$, a contradiction. So $|N_k(v_j)|\ge 1$ for any $v_j\in V_j$. Thus $G[V_j \cup V_k]$ is perfect to $V_j$ and by Corollary \ref{C43}, $|N_k(v_j)|= 1$ for any $v_j\in V_j$.

For any $v_j\in V_j\setminus N_j(V_i)$ and $v_i\in V_i$, suppose $v_j^{1}\in N_j(v_i)$ and $v\in N_k(v_j^{1})$. Then the only possible $C_4$ in $G+v_iv_j$ is $v_jv_iv_{j}^{1}vv_j$. Hence $N_k(v_j)=\{v\}$. By the arbitrariness of $v_j$ and $v_i$ we have $N_j[v]=V_j.$

	
\vspace{1em}

\noindent{\bf Proof of Theorem \ref{T15}. }By Theorem \ref{T12} and the construction of $W^{(4,3,n)}$ in Section 2, we have $sat(K_{3} ^{ n}, C_4) = 3n$. By Corollary \ref{C43} and Lemma \ref{L44}, we can assume that $G \, [ \, V_1 \cup V_2 \, ]$ is perfect to $V_1$ and non-perfect to $V_2$. By Lemma \ref{L45}, there is $z_3\in V_3$ such that $N_3(x)=\{z_3\}$ for any $x\in V_2$. Thus $G \, [ \, V_2 \cup V_3 \, ]$ is perfect to $V_2$ and non-perfect to $V_3$  by Corollary \ref{C43}. By Lemma \ref{L45}, there is $z_1\in V_1$ such that $N_1(y)=\{z_1\}$ for any $y\in V_3$. By the same argument, there is $z_2\in V_2$ such that $N_2(z)=\{z_2\}$ for any $z\in V_1$.
 Thus $G=W^{(4,3,n)}$ and $Sat(K_{3} ^{ n}, C_4) = \Omega^{(4,3,n)}$. \qed

\subsection{Non-tripartite}

We first establish the upper bound of Theorem \ref{T19} by construct a graph family $\{Y_k\}$.

	Let $Y_4$ be a $4$-partite graph with $n$ vertices in each part, $Y_4 = V_1 \cup V_2 \cup V_3 \cup V_4$. Pick $x_i \in V_i$ for $i \in [3]$ and $X=\{x_1,x_2,x_3\}$, $G[X]=K_3$. Let $B_i =V_i \setminus X$. For any $u \in B_1$, $N(u)=\{x_2\}$. For any $u \in B_4=V_4$, $N(u)=\{x_3\}$. $G[B_2 \cup B_3] $ is a perfect matching and for any $u \in B_2 \cup B_3$, $|N(u)|=2$, $N(u) \setminus \left( B_2 \cup B_3 \right) = \{x_1\}$. The image of $Y_4$ is in Figure 14. For $k \geq 5$, We would construct $Y_k$ by induction.
\vspace{0.5em}

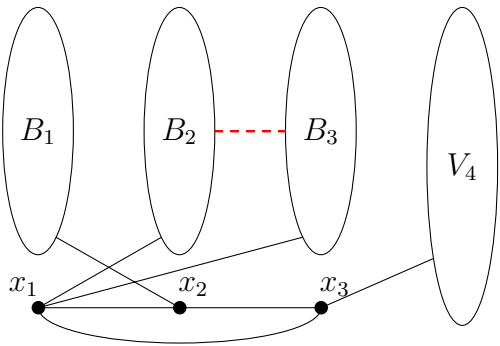
\begin{figure}[H]
\centering
\begin{tikzpicture}[scale=.47]

\draw (-10,0) arc(0:360:1cm and 3.5cm) ;
\draw (-6,0) arc(0:360:1cm and 3.5cm) ;
\draw (-2,0) arc(0:360:1cm and 3.5cm) ;
\draw (2,-1) arc(0:360:1cm and 4.5cm) ;

\draw (1,-1) node[align=center]{$V_4$};
\draw (-3,0) node[align=center]{$B_3$};
\draw (-7,0) node[align=center]{$B_2$};
\draw (-11,0) node[align=center]{$B_1$};

\filldraw (-2.99,-5) circle (5pt);
\filldraw (-6.99,-5) circle (5pt);
\filldraw (-10.99,-5) circle (5pt);

\draw (-2.6,-4.4) node[align=center]{$x_3$};
\draw (-6.6,-4.4) node[align=center]{$x_2$};
\draw (-11.4,-4.4) node[align=center]{$x_1$};


\draw (-2.99,-5)--(-6.99,-5) (-10.99,-5)--(-6.99,-5)  (-2.99,-5) arc(0:-180: 4cm and 1cm) (-10.99,-5);

\draw (-10.5,-3) --(-6.99,-5) (-7.5,-3) --(-10.99,-5) (-3.5,-3) --(-10.99,-5) (0.2,-3.6) --(-2.99,-5) ;
\draw[red,dashed,thick]   (-6,0) -- (-4,0) (-6,0) -- (-4,0) (-6,0) -- (-4,0) ;

\end{tikzpicture}\\

\caption{\centering $Y_4$. The black solid lines and arcs represent the complete connection between vertices, and the red dashed line represents the perfect matching. }

\end{figure}

\begin{definition}\label{D46}
For $k \geq 5$,  \\
\textbf{(i)} \,$Y_5=V_1 \cup V_2 \cup ...  \cup V_5$ is a $5$-partite graph with $n$ vertices in each part. $G[V_1 \cup ... \cup V_4] =Y_4$, and all the vertices in $V_5$ are connected to the only vertex $x_1 \in V_1$ where $x_1$ is the same as in $G[V_1 \cup ... \cup V_4]=Y_4$. (Figure 15) \\
\textbf{(ii)} \,$Y_6=Y_6^\delta$ comes in two different shapes depending on the parity of $n$, where $\delta \equiv n \, (\, mod \; 2 \, )$. $Y_6^\delta = V_1 \cup V_2 \cup ... \cup V_6$ is a $6$-partite graph with $n$ vertices in each part. Pick $x_i \in V_i$ for $i \in [3]$ and $X=\{x_1,x_2,x_3\}$, $G[X]=K_3$. Let $B_i =V_i \setminus X$, then $|B_i|=n-1$ for $i \in [3]$.

	(1) $\delta=1$ such that $n$ is odd. Let $i,j,k $ be distinct elements of $ \{1,2,3\}$. Partition $B_i$ into $B_i^j,B_i^k $ such that $|B_i^j|=|B_i^k|=\frac{n-1}{2}$. $G[B_i^k \cup B_j^k]$ is a perfect matching and for any $u \in B_i^k \cup B_j^k$, $|N(u)|=2$, $N(u) \setminus \left( B_i^k \cup B_j^k \right) = \{x_k\}$. For any $u \in V_4, \, N(u)=\{x_2\}$. For any $u \in V_5, \, N(u)=\{x_3\}$. For any $u \in V_6, \, N(u)=\{x_1\}$. (Figure 16)

	(2) $\delta=2$ such that $n$ is even. For $i \in \{2,3\}$, partition $B_i$ into $B_i^+,B_i^- $ such that $|B_i^+|=\frac{n}{2}$ and $|B_i^-|=\frac{n}{2}-1$. Partition $V_4$ into $V_4^2,V_4^3 $ such that $|V_4^2|=|V_4^3|=\frac{n}{2}$. For $i \in \{2,3\}$, $G[B_i^+ \cup V_4^i]$ is a perfect matching. $G[B_2^- \cup B_3^-]$ is a perfect matching. For any $u \in B_2 \cup B_3 \cup V_4 $, $d(u)=2$ and $N(u) \setminus \left( B_2 \cup B_3 \cup V_4 \right) =\{x_1\}$. For any $u \in B_1, \, N(u)=\{x_2\}$. For any $u \in V_5, \, N(u)=\{x_3\}$. For any $u \in V_6, \, N(u)=\{x_1\}$. (Figure 17)
\vspace{0.2em}

\noindent\textbf{(iii)} For $k \geq 7$, $Y_k=Y_{k-2} \cup V_{k-1} \cup V_k$ is a $k$-partite graph with $n$ vertices in each part. $G[V_{k-1} \cup V_k]$ is a perfect matching, and all vertices in $V_{k-1} \cup V_k$ are connected to $x_1 \in V_1$ where $x_1$ is the same as in $Y_{k-2}$. $d(u)=2$ for all $u \in V_{k-1} \cup V_k $. (Figure 18)
\qed
\end{definition}	

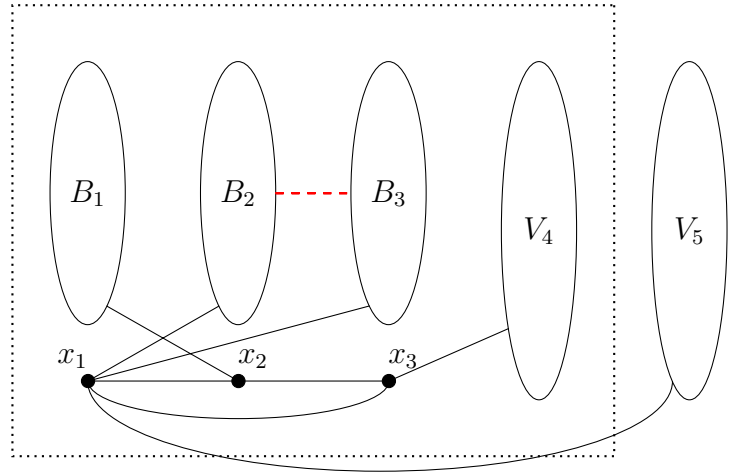
\begin{figure}[H]
\centering
\begin{tikzpicture}[scale=.5]

\draw (-10,0) arc(0:360:1cm and 3.5cm) ;
\draw (-6,0) arc(0:360:1cm and 3.5cm) ;
\draw (-2,0) arc(0:360:1cm and 3.5cm) ;
\draw (2,-1) arc(0:360:1cm and 4.5cm) ;
\draw (6,-1) arc(0:360:1cm and 4.5cm) ;

\draw (5,-1) node[align=center]{$V_5$};
\draw (1,-1) node[align=center]{$V_4$};
\draw (-3,0) node[align=center]{$B_3$};
\draw (-7,0) node[align=center]{$B_2$};
\draw (-11,0) node[align=center]{$B_1$};

\filldraw (-2.99,-5) circle (5pt);
\filldraw (-6.99,-5) circle (5pt);
\filldraw (-10.99,-5) circle (5pt);

\draw (-2.6,-4.4) node[align=center]{$x_3$};
\draw (-6.6,-4.4) node[align=center]{$x_2$};
\draw (-11.4,-4.4) node[align=center]{$x_1$};


\draw (-2.99,-5)--(-6.99,-5) (-10.99,-5)--(-6.99,-5)  (-2.99,-5) arc(0:-180: 4cm and 1cm) (-10.99,-5);

\draw (-10.5,-3) --(-6.99,-5) (-7.5,-3) --(-10.99,-5) (-3.5,-3) --(-10.99,-5) (0.2,-3.6) --(-2.99,-5) (4.55,-5) arc(0:-180: 7.775cm and 2.4cm) (-14.99,-5) ;
\draw[red,dashed,thick]   (-6,0) -- (-4,0) (-6,0) -- (-4,0) (-6,0) -- (-4,0) ;

\draw[dotted,thick] (3,5) -- (3,-7) -- (-13,-7) -- (-13,5) -- (3,5);

\end{tikzpicture}\\

\caption{\centering $Y_5$. The black solid lines and arcs represent the complete connection between vertices, and the dotted rectangle represents $Y_4$. }

\end{figure}

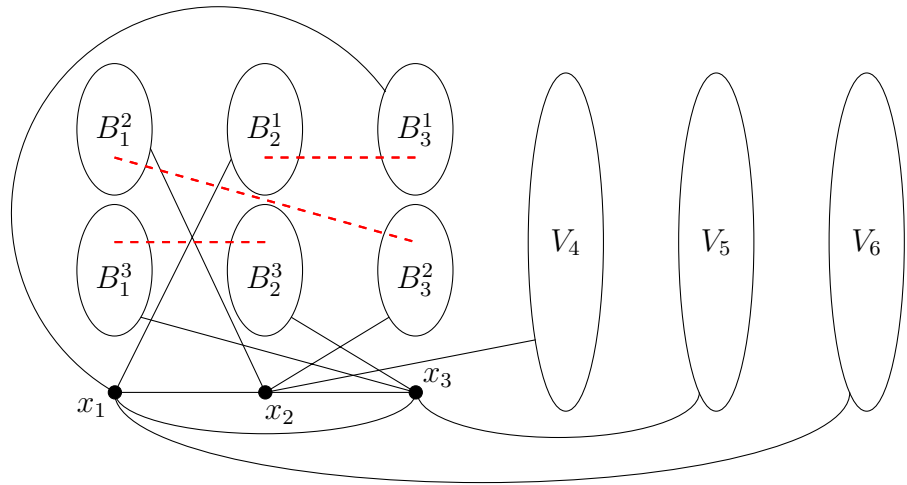
\begin{figure}[H]
\centering
\begin{tikzpicture}[scale=.5]

\draw (-10,-1.75) arc(0:360:1cm and 1.75cm) ;
\draw (-6,-1.75) arc(0:360:1cm and 1.75cm) ;
\draw (-2,-1.75) arc(0:360:1cm and 1.75cm) ;
\draw (-10,2) arc(0:360:1cm and 1.75cm) ;
\draw (-6,2) arc(0:360:1cm and 1.75cm) ;
\draw (-2,2) arc(0:360:1cm and 1.75cm) ;
\draw (2,-1) arc(0:360:1cm and 4.5cm) ;
\draw (6,-1) arc(0:360:1cm and 4.5cm) ;
\draw (10,-1) arc(0:360:1cm and 4.5cm) ;


\draw (9,-1) node[align=center]{$V_6$};
\draw (5,-1) node[align=center]{$V_5$};
\draw (1,-1) node[align=center]{$V_4$};
\draw (-3,2) node[align=center]{$B_3^1$};
\draw (-7,2) node[align=center]{$B_2^1$};
\draw (-11,2) node[align=center]{$B_1^2$};
\draw (-3,-2) node[align=center]{$B_3^2$};
\draw (-7,-2) node[align=center]{$B_2^3$};
\draw (-11,-2) node[align=center]{$B_1^3$};

\filldraw (-2.99,-5) circle (5pt);
\filldraw (-6.99,-5) circle (5pt);
\filldraw (-10.99,-5) circle (5pt);

\draw (-2.4,-4.6) node[align=center]{$x_3$};
\draw (-6.6,-5.6) node[align=center]{$x_2$};
\draw (-11.6,-5.4) node[align=center]{$x_1$};


\draw (-2.99,-5)--(-6.99,-5) (-10.99,-5)--(-6.99,-5)  (-2.99,-5) arc(0:-180: 4cm and 1.1cm) (-10.99,-5);

\draw  (-10.3,-3) -- (-2.99,-5)   (-6.3,-3) -- (-2.99,-5) (-3.7,-3) --(-6.99,-5);
\draw  (-10.05,1.5) -- (-6.99,-5)  (-10.99,-5) -- (-7.9,1.2);
\draw  (-10.99,-5) arc(240:36: 5.5cm and 5.5cm);

\draw (0.2,-3.6) --(-6.99,-5) (4.55,-5) arc(0:-180: 3.75cm and 1.2cm) (-2.99,-5) (8.55,-5) arc(0:-180: 9.775cm and 2.4cm) (-14.99,-5);

\draw[red,dashed,thick] (-11,-1) -- (-6.8,-1)  (-11,1.25) -- (-3,-1) (-3,1.25) -- (-7.2,1.25);
\draw[red,dashed,thick] (-11,-1) -- (-6.8,-1)  (-11,1.25) -- (-3,-1) (-3,1.25) -- (-7.2,1.25);


\end{tikzpicture}\\

\caption{\centering $Y_6^1$. The black solid lines and arcs represent the complete connection between vertices.  The red dashed lines represent the perfect matching.}
\end{figure}

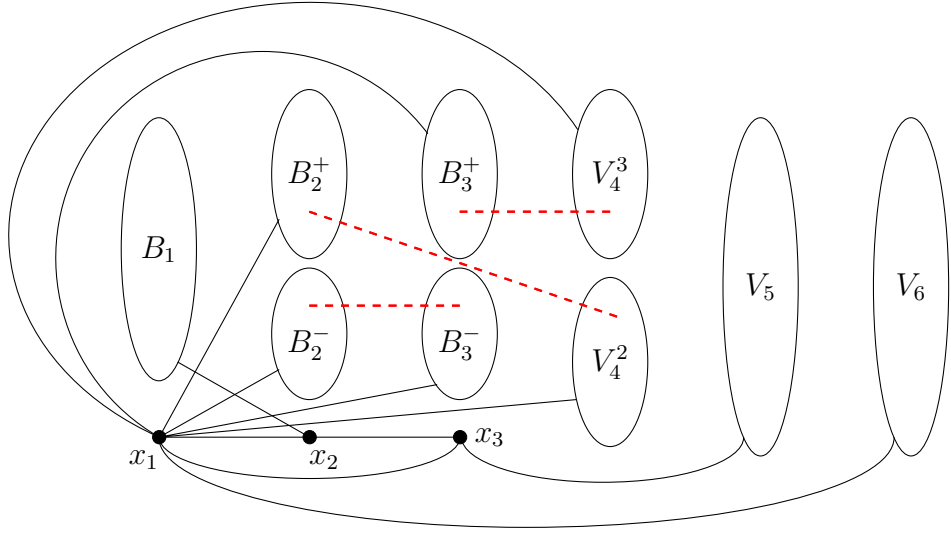
\begin{figure}[H]
\centering
\begin{tikzpicture}[scale=.5]

\draw (-10,0) arc(0:360:1cm and 3.5cm) ;
\draw (-6,2) arc(0:360:1cm and 2.25cm) ;
\draw (-2,2) arc(0:360:1cm and 2.25cm) ;
\draw (-6,-2.25) arc(0:360:1cm and 1.75cm) ;
\draw (-2,-2.25) arc(0:360:1cm and 1.75cm) ;
\draw (2,2) arc(0:360:1cm and 2.25cm) ;
\draw (2,-3) arc(0:360:1cm and 2.25cm) ;
\draw (6,-1) arc(0:360:1cm and 4.5cm) ;
\draw (10,-1) arc(0:360:1cm and 4.5cm) ;

\draw (9,-1) node[align=center]{$V_6$};
\draw (5,-1) node[align=center]{$V_5$};
\draw (1,2) node[align=center]{$V_4^3$};
\draw (1,-3) node[align=center]{$V_4^2$};
\draw (-3,2) node[align=center]{$B_3^+$};
\draw (-7,2) node[align=center]{$B_2^+$};
\draw (-3,-2.5) node[align=center]{$B_3^-$};
\draw (-7,-2.5) node[align=center]{$B_2^-$};
\draw (-11,0) node[align=center]{$B_1$};

\filldraw (-2.99,-5) circle (5pt);
\filldraw (-6.99,-5) circle (5pt);
\filldraw (-10.99,-5) circle (5pt);

\draw (-2.2,-5) node[align=center]{$x_3$};
\draw (-6.6,-5.6) node[align=center]{$x_2$};
\draw (-11.4,-5.6) node[align=center]{$x_1$};


\draw (-2.99,-5)--(-6.99,-5) (-10.99,-5)--(-6.99,-5)  (-2.99,-5) arc(0:-180: 4cm and 1.1cm) (-10.99,-5) ;

\draw (-10.5,-3) --(-6.99,-5)  (4.55,-5) arc(0:-180: 3.75cm and 1.2cm)  (8.55,-5) arc(0:-180: 9.775cm and 2.4cm)  ;

\draw  (-7.8,-3.2) --(-10.99,-5) (-3.6,-3.6) --(-10.99,-5) (0.1,-4) --(-10.99,-5) --  (-7.8,0.8);

\draw (-10.99,-5) arc(240:37: 5.5cm and 5.5cm);
\draw (-10.99,-5) arc(240:27: 8cm and 6.2cm);

\draw[red,dashed,thick] (-7,-1.5) -- (-3,-1.5)  (-7,1) -- (1.2,-1.8) (1,1) -- (-3,1);
\draw[red,dashed,thick] (-7,-1.5) -- (-3,-1.5)  (-7,1) -- (1.2,-1.8) (1,1) -- (-3,1);


\end{tikzpicture}\\

\caption{\centering $Y_6^2$. The black solid lines and arcs represent the complete connection between vertices.  The red dashed lines represent the perfect matching. }
\end{figure}

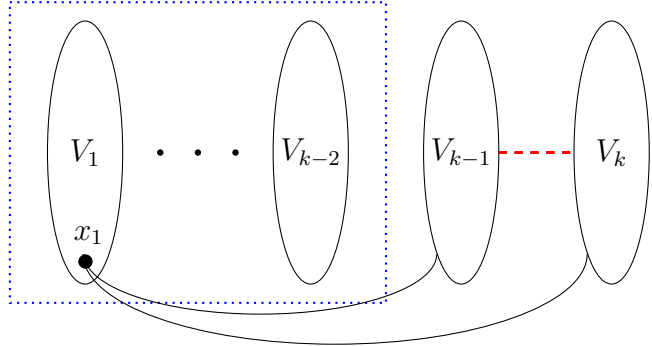
\begin{figure}[H]
\centering

\begin{tikzpicture}[scale=.5]

\draw (-16,0) arc(0:360:1cm and 3.5cm) ;
\draw (-10,0) arc(0:360:1cm and 3.5cm) ;
\draw (-6,0) arc(0:360:1cm and 3.5cm) ;
\draw (-2,0) arc(0:360:1cm and 3.5cm) ;

\draw (-3,0) node[align=center]{$V_{k}$};
\draw (-7,0) node[align=center]{$V_{k-1}$};
\draw (-11,0) node[align=center]{$V_{k-2}$};
\draw (-17,0) node[align=center]{$V_{1}$};

\filldraw (-16.99,-2.9) circle (5pt);
\filldraw (-13,0) circle (2pt) (-14,0) circle (2pt) (-15,0) circle (2pt);

\draw (-16.9,-2.2) node[align=center]{$x_1$};

\draw[red,dashed,thick] (-6,0) -- (-4,0)  (-6,0) -- (-4,0) (-6,0) -- (-4,0);

\draw[blue,dotted, thick] (-19,4) -- (-9,4) -- (-9,-4) --  (-19,-4) -- (-19,4);
\draw (-7.65,-2.7) arc(0:-180: 4.69cm and 1.6cm)  ;
\draw (-3.65,-2.7) arc(0:-180: 6.7cm and 2.4cm)  ;

\end{tikzpicture}\\

\caption{\centering $Y_k$. The black solid arcs represent the complete connection between vertices. The red dashed line represents the perfect matching. The blue dotted rectangle represents $Y_{k-2}$. }
\end{figure}

\begin{prop}\label{P47}
$Y_k$ is $C_4$-saturated. 
\end{prop}
	
\noindent{\bf Proof. }The vertices of $Y_k$ can be divided into $A,B,X$: the vertex in $A$ has degree $1$; the vertex in $B$ has degree $2$\,; $X=\{x_1,x_2,x_3\}$. It east to check that $Y_k$ is $C_4$-free. 

	The admissible nonedge $uv$ in $Y_k$ can only be\,: (1) $u,v \in A$\,; (2) $u \in A \cup B$, $v \in B$\,; (3) $v=x_j$ for some $j \in [3]$. 
	
	(1) If $u,v \in A$ then $N(u) \neq N(v)$. Let $x_i$ be the neighbour of $u$ and $x_j$ be the neighbour of $v$, the addition of $uv$ would form a $C_4 = u x_i x_j vu$. (2) If $u \in A \cup B$, $v \in B$, let $x_i \in N(u)$ and $N(v)=\{x_j,w\}$ where $w \in B$. If $i=j$, then the addition of $uv$ would form a $C_4=ux_iwvu$\,; else if $i \neq j$, the addition of $uv$ would form a $C_4=ux_ix_jvu$. (3) If $v=x_j$ for some $j \in [3]$, let $x_i \in N(u)$ and $x_w$ be the distinct elements of $x_i,x_j$ in $X$. The addition of $uv$ would form a $C_4 = u x_i x_w vu$. Hence $Y_k$ is $C_4$-saturated. \qed
\vspace{0.3em}

	Proposition \ref{P47} shows the upper bound of Theorem \ref{T19}. Clearly $e(Y_4)= 5n-1$, and for $k \geq 5$ we have (One can see that $|A|=3n$ if $k$ is even and $n$ is odd, or $|A|=3n-1$)
\begin{equation*}
\begin{aligned}
  	e(Y_k) & = e(A,X) + e(G[B]) + e(B,X) + e(G[X]) \\
  		& = |A|+\frac{|B|}{2} + |B| + |X|  = kn + \frac{|B|}{2} \\
  		& = kn + \lfloor \frac{kn-(3n-1)-3}{2} \rfloor \\
  		& = \lfloor  \frac{3(k-1)n -2}{2} \rfloor .
\end{aligned} 
\end{equation*}
Thus we have proved the upper bounds of Theorem \ref{T19}.  
\vspace{0.3em}

	Next we focus on the lower bound of $sat(K_k^n, C_4)$ when $k \geq 4$. Let $G$ be a $C_4$-saturated $k$-partite graph with $n$ vertices in each part. If $\delta(G) \geq 3$ then $e(G) \geq \frac{3}{2}kn$ and we are done, so we assume $\delta(G) \leq 2$.

\begin{lemma}\label{L48}  Let $u,v$ be the vertices of degree 1 in $G$. If $u,v$ lie in different parts of $G$, then $N(u) \neq N(v)$. Moreover, the number of parts containing a degree 1 vertex is at most 3.
\end{lemma}
\noindent{\bf Proof. }Let $u,v$ be the vertices of degree 1 lie in different parts of $G$. If $N(u)=N(v)$, then the addition of $uv$ would only create $C_3$ in $G$, a contradiction.

	If there exists $x_i \in V_i$ such that $d(x_i)=1$ with $i \in [4]$, let $N(x_i)=\{y_i\}$. Since the addition of $x_ix_j$ would create a $C_4$, $y_iy_j \in E(G)$ and hence $G[y_1,y_2,y_3,y_4]$ is a $K_4$, a contradiction. \qed
\vspace{0.3em}

	By Lemma \ref{L48} we know that $G$ contains at most $3n$ vertices of degree 1. Moreover, let $x_1,x_2 \in V_i$ for some $i$ such that $N(x_1)=\{y_1\}, N(x_2)=\{y_2\}$ with $y_1 \neq y_2$. Then a new graph $G'= \left( G - x_2y_2 \right) + x_2y_1$ is still a $C_4$-saturated $k$-partite graph with the same number of edges as $G$. So we assume that all degree 1 vertices in the same part of $G$ have the same neighbour.
\vspace{0.3em}

\begin{prop}\label{P49}
If $\delta(G)=1$, then $e(G) \geq \lfloor  \frac{3(k-1)n -2}{2} \rfloor$.
\end{prop}
	
\noindent{\bf Proof. }Let $x \in V_1$ be the degree 1 vertex and $N(x)=\{y\}$. Let $S_i = \{ u \in V(G): d(x,u)=i \, \}$. Since $diam_p(G) \leq 3$, we only need to consider $S_1=\{y\},S_2,S_3$ and possible $S_4 \subseteq V_1$. By our assumption all degree 1 vertices in $V_1$ would connected to $y$, let $A_2$ be the set of such vertices in $S_2$ and $B_2=S_2 \setminus A_2$. Then for any $u \in S_4$ we have $d(u) \geq 2$. By Lemma \ref{L48} the degree 1 vertices not in $V_1$ are all belong to $S_3$. Let $A_3=\{u \in S_3 : d(u)=1\}$ and $B_3=S_3 \setminus A_3$. By Lemma \ref{L48} we have $|A_3| \leq 2n$, and moreover $|S_4|+|A_2| \leq |V_1|-1 = n-1$.

	For a vertex $u \neq x,y$, let $M_i(u)=N(u) \cap S_i$ for $i \in [4]$. Clearly if $u \in S_j$ for some $j$, then if $|i-j| \geq 2$ we have $M_i(u)=\emptyset$. Define a weight function $f$ as\,:
\begin{equation*}
f(u)=\left\{
\begin{aligned}
 \; & |M_1(u)| + \frac{1}{2} |M_2(u)| \, , \; & u \in S_2 \, ,\\
 \; & |M_2(u)| + \frac{1}{2} |M_3(u)| + \frac{1}{2}| M_4(u)| \, , \; & u \in S_3 \, , \\\
 \; & \frac{1}{2} |M_3(u)| \, , \; & u \in S_4 \, .
\end{aligned}
\right.
\end{equation*}
Then $\sum_{u \neq x,y} f(u) = e(G)-1$. 
\vspace{0.5em}

\noindent {\bf Claim 1. }For any $u \in B_2 \cup B_3$, $f(u) \geq \frac{3}{2}$.

	If $u \in B_2$, since $ux$ is connected by a $P_4$, $u$ must have at least one neighbour in $S_2$ and hence $|M_1(u)|+\frac{1}{2}|M_2(u)| \geq \frac{3}{2}$. If $u \in B_3$, since $d(u) \geq 2$, either $|M_2(u)| \geq 2$ or $|M_2(u)|=1 $, $|M_3(u)| + |M_4(u)| \geq 1$, so $f(u) \geq \frac{3}{2}$. \q
\vspace{0.5em}

\noindent {\bf Claim 2. }For any $u \in A_2 \cup A_3 \cup S_4$, $f(u) \geq 1$. 

	If $u \in A_2 \cup A_3$, by our definition we immediately have the result. If $u \in S_4$, since $d(u) \geq 2$ and $S_4 \subseteq V_1$, all neighbour of $u$ must lie in $S_3$, so $f(u) \geq \frac{1}{2}|M_3(u)| \geq 1$. \q
\vspace{0.5em}

	By Claims 1 and 2 we can calculate the edges in $G$. Recall that $|A_3| \leq 2n$, $|S_4|+|A_2| \leq   n-1$ and $|A_2|+ |B_2|+|A_3| +|B_3|+|S_4|=kn-2$, we have 
\begin{equation*}
\begin{aligned}
e(G) - 1 & = \sum_{u \neq x,y} f(u) \\
& = \sum_{u \in B_2 \cup B_3} f(u) + \sum_{u \in A_2 \cup A_3 \cup S_4} f(u) \\
& \geq \frac{3}{2} (|B_2|+|B_3|) + (|A_2|+|A_3|+|S_4|) \\
& = \frac{3}{2}(kn-2) - \frac{1}{2} (|A_2|+|A_3|+|S_4|)  \\
& \geq \frac{3}{2}(kn-2) - \frac{1}{2}(3n-1) = \frac{3(k-1)n-5}{2} \, ,
\end{aligned}
\end{equation*}
so we have $e(G) \geq \frac{3(k-1)n-3}{2}$. Since $e(G)$ is an integer we have $e(G) \geq \lfloor  \frac{3(k-1)n -2}{2} \rfloor$ and we are done. \qed
\vspace{0.5em}

\begin{prop}\label{P410}
If $\delta(G)=2$ and there exists a degree 2 vertex not in a triangle, then we have $e(G) \geq \frac{(3k-1)n-4}{2}$.
\end{prop}

\noindent{\bf Proof. }Let $x \in V_1$ be the degree 2 vertex and $N(x)=\{y_1,y_2\}$ such that $y_1y_2 \not \in E(G)$. Similar to the proof of Proposition \ref{P49}, let $S_i=\{ u \in V(G) : d(x,u)=i\}$. Since $diam_p(G) \leq 3$, we only consider $S_1=\{y_1,y_2\},S_2,S_3$ and possible $S_4 \subseteq V_1$. We have $|S_2|+|S_3|+|S_4|=kn-3$ and $|S_4| \leq n-1$.
	
	For a vertex $u \neq x,y_1,y_2$, let $M_i(u)=N(u) \cap S_i$ for $i \in [4]$. Define the same weight function $f$ as\,:
\begin{equation*}
f(u)=\left\{
\begin{aligned}
 \; & |M_1(u)| + \frac{1}{2} |M_2(u)| \, , \; & u \in S_2 \, ,\\
 \; & |M_2(u)| + \frac{1}{2} |M_3(u)| + \frac{1}{2}| M_4(u)| \, , \; & u \in S_3 \, , \\\
 \; & \frac{1}{2} |M_3(u)| \, , \; & u \in S_4 \, .
\end{aligned}
\right.
\end{equation*}
Then $\sum_{u \neq x,y_1,y_2} f(u) = e(G)-2$. 
\vspace{0.5em}

\noindent {\bf Claim. }For any $u \in S_2 \cup S_3$, $f(u) \geq \frac{3}{2}$.

	If $u \in B_2$, since $ux$ is connected by a $P_4$ and $y_1y_2$ is not an edge, $u$ must have at least one neighbour in $S_2$ and hence $|M_1(u)|+\frac{1}{2}|M_2(u)| \geq \frac{3}{2}$. If $u \in B_3$, since $d(u) \geq 2$, either $|M_2(u)| \geq 2$ or $|M_2(u)|=1 $, $|M_3(u)| + |M_4(u)| \geq 1$, so $f(u) \geq \frac{3}{2}$. \q
\vspace{0.5em}

	By a similar argument we can prove that for any $u \in S_4$ we have $f(u) \geq 1$. Thus the number of edges in $G$ is at least

\begin{equation*}
\begin{aligned}
e(G)  & = \sum_{u \neq x,y_1,y_2} f(u) + 2 \\
& = \sum_{u \in S_2 \cup S_3} f(u) + \sum_{u \in S_4 } f(u) +2  \\
& \geq \frac{3}{2} (|S_2|+|S_3|) + |S_4| +2 \\
& = \frac{3}{2}(kn-3) - \frac{1}{2} |S_4| + 2 \\
& \geq \frac{3}{2}(kn-3) - \frac{1}{2}(n-1)+2 = \frac{(3k-1)n-4}{2} \, .
\end{aligned}
\end{equation*}
\qed

	The next Lemma is a direct result in \cite{FI}. Although they only proved the case in complete graph $K_{kn}$, such result also holds for complete multipartite graph $K_k^n$. We would restate their proof.

\begin{lemma}\label{L411}\text{\bf \cite{FI}}\;
Let $H$ be a connected multipartite graph. If each edge of $H$ is in a triangle, then $e(H) \geq \frac{3}{2}(|H|-1)$.
\end{lemma}
\noindent {\bf Proof. }Pick a spanning tree $T$ with $|H|-1$ edges. An edge of $H-T$ can form a triangle with at most two edges of $T$, so $H-T$ has at least $\frac{1}{2}(|H|-1)$ edges. \qed
\vspace{0.8em}

\begin{prop}\label{P412}
If $\delta(G)=2$ and each degree 2 vertex is in a triangle, then we have $e(G) \geq \frac{3(k-1)n}{2}$.
\end{prop}
\noindent {\bf Proof. }Let $H$ be the subgraph of all $C_3$'s containing at least one degree 2 vertex. Let $H_1,...,H_m$ be the components of $H$. By Lemma \ref{L411} we have $e(H_i) \geq \frac{3}{2}(|H_i|-1)$. For distinct $H_i,H_j$, if there exists degree 2 vertices $u \in V(H_i)$ and $v \in V(H_j)$ such that $u,v$ lie in different parts, then there exists at least an edge between $H_i,H_j$ since $u,v$ is connected by a $P_4$. 

	For any $i\in[k]$, suppose that there exists $m_i$ components in $H$ such that these $m_i$ components only contain degree 2 vertices in $V_i$. Clearly $m_1+...+m_k \leq m$, $|m_i| \leq n$ and there are at most $M$ pairs of $(H_i,H_j)$ which are not connected by an edge, where $M=\sum_{1 \leq i \leq m} \binom{m_i}{2} \leq \binom{m}{2}$. Hence
\begin{equation*}
e(G[V(H)]) \geq \sum_{i=1}^m \frac{3}{2}(|H_i|-1) + \binom{m}{2} - \sum_{1 \leq i \leq m} \binom{m_i}{2}  =\,: \mathscr{S}. 
\end{equation*}

	If $m \leq 2n$, then $\mathscr{S} \geq \sum \limits_{i=1}^m \frac{3}{2}(|H_i|-1) \geq  \frac{3}{2}|H|-3n$. If $m>2n$, assume that $rn<m \leq (r+1)n$ for some integer $r$. Since $m_1+...+m_k \leq m$, $|m_i| \leq n$ and the combinatorial binomial is a convex function, we have 
\begin{equation*}
\begin{aligned}
\binom{m}{2} - \sum_{1 \leq i \leq m} \binom{m_i}{2} -m & \geq \frac{m(m-1)}{2} - (r+1)\frac{n(n-1)}{2} - (r+1)n \\
& > \frac{(r^2 n^2-rn) - (r+1)n(n+1)}{2} = \frac{(r^2-r-1) n^2-(2r+1)n}{2} \\
& \geq \frac{(r^2-r-1) n-(2r+1)n}{2} = \frac{(r^2-3r-2) n}{2} > -3n \, ,
\end{aligned}
\end{equation*}
hence we have $\mathscr{S} = \sum_{i=1}^m \frac{3}{2}|H_i|- m + \binom{m}{2} - \sum_{1 \leq i \leq m} \binom{m_i}{2} \geq  \frac{3}{2}|H|-3n$.
	
	
	Since all degree 2 vertices are in $H$, then for any $u \in V(G) \setminus V(H)$ we have $d(u) \geq 3$. So we have $e(G) \geq e(G[V(H)])+ \frac{3}{2}(kn-|H|) = \frac{3(k-1)n}{2}$. \qed

\vspace{0.8em}

\noindent{\bf Proof of Theorem \ref{T19}. }By Propositions \ref{P47}, \ref{P49}, \ref{P410} and \ref{P412} we are done. \qed

\section{$C_5$-saturation multipartite graph}

 Recall $diam_p(G)$ is the maximum distance between two vertices not lie in the same part of $G$.

Now we construct a sequence $G=G_0\supseteq G_1\supseteq\ldots$ of induced subgraphs of $G$ as follows. If $G_i$ has a vertex $v_i$ of degree 1, we let $G_{i+1}:=G_i-v_i$; if not, we terminate our sequence and set $Tr(G)=G_i$. Let $Br(G)=G[V(G) \setminus V(Tr(G)]$. We call $Tr(G)$ and $Br(G)$ the  \textit{trunk} and \textit{branch} of $G$ respectively. Without confusion, we abbreviate $Tr(G)$ and $Br(G)$ as $Tr$ and $Br$. We now have the following  observations.
\vspace{0.1em}

\begin{observation}\label{O51} Let $G$ be a connected graph but  not a tree. Assume that $Br_1,\ldots,Br_m$ are components of $Br$. Then \\
\indent \text{(i)} $V(Tr) \neq \emptyset$ and $Tr$ is connected.\\
\indent \text{(ii)} $\delta(Tr) \geq 2$. Moreover, if $u$ lies in a cycle of $G$, then $u \in V(Tr)$. \\
\indent \text{(iii)} $e(Tr)=e(G)-|Br|$. \\
\indent \text{(iv)}  For  $1 \leq j \leq m$, $Br_j$ is a tree, $N(V(Br_j)) \subseteq V(Tr)$, $|N(V(Br_j))| =1 $, and $\overline{Br}_j: = G \left[ N[V(Br_j)] \right]$ is a tree.
\end{observation}

	The first three properties are easy to check by the construction. Still it is easy to see that $Br_j$ is a tree, $N(V(Br_j)) \subseteq V(Tr)$ and $|N(V(Br_j))| =1 $.  If $\overline{Br}_j$ is not a tree, then some vertex of $Br_j$ lies in a cycle of $G$, a contradiction with (ii).
	
	By Observation \ref{O51} (iv), we let $N(V(Br_j))=\{\overline{r}_j\}$, $1\le j\le m$. For $1\le j\le m$, the distance $d(u,Tr)$ of a vertex $u \in V(Br_j)$ to $Tr$ is the number of edges in the only path  connecting $u $ and $\overline{r}_j$, and $rad(Br_j) = \max_{u \in V(Br_j)} d(u,Tr)$ is the \textit{radius} of $Br_j$.
	
	We now focus on the lower bound of $sat(K_4^n, C_5)$. Let $G \in Sat( K_4^n, C_5)$.
Assume $Tr$ is the trunk of $G$ and $Br$ is the branch of $G$ with $m$  components $Br_1,\ldots,Br_m$. By Theorem \ref{T12}, $e(G) \geq 4n=|G|$. Note that $G$ is connected, then $G$ has at least one cycle. By Observation \ref{O51} (i), $Tr$ is connected. Since $G$ is $C_5$-saturated,  $diam_p(G) \leq 4$ and $Tr $ is $C_5$-saturated. We have following properties.
\vspace{0.1em} 	
	
\begin{lemma}\label{L52}  For any $j \in [m]$, $rad(Br_j) \leq 2$.
\end{lemma}

\noindent{\bf Proof. }Suppose  that there exists $Br_j$ with $rad(Br_j) \geq 3$. Then there exists a vertex $x \in V(Br_j)$ such that $xyz\overline{r}_j$ is a path in $\overline{Br}_j$ with $d(x,\overline{r}_j)=3$. We have either $xz$ or $x\overline{r}_j$ is an admissible nonedge.  By Observation \ref{O51} (iv), $G+xz$ or $G+x\overline{r}_j$ would only create a $C_3$ or a $C_4$, a contradiction.\qed
\vspace{0.5em}

By the proof of Lemma \ref{L52}, we have the following corollary.

\begin{cor}\label{C53}  If $rad(Br_j) = 2$, then all the vertices $u \in V(Br_j)$ with $d(u,Tr)=2$  lie in the same part as $\overline{r}_j$.
\end{cor}
\vspace{0.0em}

\begin{lemma}\label{L54}  Let $u \in V(Tr)$. Then $|\{j|rad(Br_j) = 2,~\overline{r}_j=u,~1\le j\le m\}|\leq 1$.
\end{lemma}

\noindent{\bf Proof. }Suppose there are $j_1,j_2\in\{j|rad(Br_j) = 2,~\overline{r}_j=u,~1\le j\le m\}$. Let $x_iy_i u$ be the path of $\overline{Br}_{j_i}$, where $i=1,2$. By Corollary \ref{C53} $x_1,x_2,u$ lie in the same part. Thus $x_1y_2$ is an admissible nonedge of $G$. It is easy to see that $G[V(Br_{j_1}) \cup V(Br_{j_2}) \cup \{u\}]$ is a tree. Then  $G+x_1y_2$ can only have a $C_4=x_1y_1uy_2x_1$, a contradiction. \qed
\vspace{0.0em}

\begin{lemma}\label{L55}  $|\{j|rad(Br_j) = 2,~1\le j\le m\}|\le 1$.
\end{lemma}

\noindent{\bf Proof. }Suppose there are $j_1,j_2\in\{j|rad(Br_j) = 2,~1\le j\le m\}$. By Lemma \ref{L54}, $\overline{r}_{j_1} \neq \overline{r}_{j_2}$. Let  $x_iy_i \overline{r}_{j_i}$ be the path of $\overline{Br}_{j_i}$, $i=1,2$. If $\overline{r}_{j_1}$ and $\overline{r}_{j_2}$ lie in the same part of $G$, say $V_1$, then $d(\overline{r}_{j_1},\overline{r}_{j_2})\ge 2$. By Corollary \ref{C53}, $x_1,x_2 \in V_1$ and $y_2 \not \in V_1$. Then  $d(x_1,y_2)\ge 5$, a contradiction to $diam_p(G) \leq 4$. Now assume $\overline{r}_{j_1}\in V_1$ and $\overline{r}_{j_2}\in V_2$. By Corollary \ref{C53}, $x_i \in V_i$ for $i=1,2$. Then  $d(x_1,x_2)\ge 5$, a contradiction to $diam_p(G) \leq 4$ again.
\qed

\begin{lemma}\label{L56}  Let $C_p$ be a cycle of $G$ such that all edges in this cycle are not belong to any other cycles, then $p \leq 7$.
\end{lemma}

\noindent{\bf Proof. }Suppose not, let $C_p=x_1x_2 \ldots x_px_1$ ($p \geq 8$) be a cycle in $G$ satisfying the condition. Then either $x_1x_3$ or $x_1x_4$ is an admissible nonedge. But $G+x_1x_3$ or $G+x_1x_4$ has no $C_5$, a contradiction.
\qed
\vspace{0.5em}

By Observation \ref{O51} (i) and (ii),  we easily have the following result.

\begin{lemma}\label{L57}
 If \,$ |Tr|\leq e(Tr) \leq |Tr|+1$, then $Tr$ can only be one of the following four types of graphs (see Figure 19): (1) a cycle; (2) a $\theta$-graph; (3) a dumb; (4) a bow.
\end{lemma}

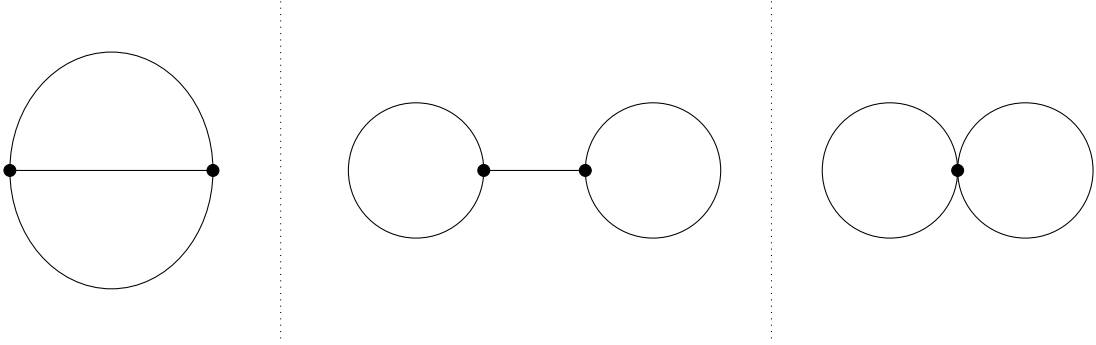
\begin{figure}[H]
\centering
\begin{tikzpicture}[scale=.45]

\draw (-10,1) arc(0:360:3cm and 3.5cm) ;
\draw (-16,1) -- (-10,1) ;
\filldraw (-16,1) circle (5pt);
\filldraw (-10,1) circle (5pt);

\draw[dotted] (-8,6) -- (-8,-4);

\draw (-2,1) arc(0:360:2cm and 2cm) ;
\draw (5,1) arc(0:360:2cm and 2cm) ;
\draw (-2,1) -- (1,1) ;
\filldraw (-2,1) circle (5pt);
\filldraw (1,1) circle (5pt);

\draw[dotted] (6.5,6) -- (6.5,-4);

\draw (12,1) arc(0:360:2cm and 2cm) ;
\draw (16,1) arc(0:360:2cm and 2cm) ;
\filldraw (12,1) circle (5pt);

\end{tikzpicture}\\

\caption{\centering a $\theta$-graph (left), a dumb (middle), a bow (right). A $\theta$-graph is a cycle with an extra path connecting two vertices on the cycle. A dumb is a graph with two cycles connecting by a path. A bow is a graph of two cycles intersecting at exactly one vertex. }

\end{figure}

	Now we are able to prove Theorem \ref{T17}.
\vspace{0.5em}	

\noindent{\bf Proof of Theorem \ref{T17}. }Let $G \in Sat(  K_{4 }^{ n}, C_5)$ and $G=V_1\cup V_2\cup V_3\cup V_4$. Assume $Tr$ is the trunk of $G$ and $Br$ is the branch of $G$ with $m$  components $Br_1,\ldots,Br_m$ and $N(V(Br_j))=\{\overline{r}_j\}$ for $1\le i\le m$.  By Theorem \ref{T12}, $e(G) \geq 4n=|G|$.           It is easy to check that $W_{\pi_i}^{(5,4,n)}$ (Figure 3) is $C_5$-saturated 4-partite graph. Thus $e(G)=sat(K_4^n, C_5) \leq e (W_{\pi_i}^{(5,4)}) = 4n+2$. Suppose  $e(G) < 4n+2$. Then $4n \leq e(G) \leq 4n+1$. By Observation \ref{O51} (iii), $|Br|=|G|-|Tr|$ and $|G|=4n$, we have $ |Tr|\leq e(Tr) \leq |Tr|+1$.
\vspace{0.5em}	

\noindent{\bf Claim 1. }$|Tr| \leq 15$.
\vspace{0.0em}

\noindent{\bf Proof of Claim 1. }Note that $Tr $ is $C_5$-saturated. By Lemma \ref{L56}, if $Tr$ is a cycle (resp. a bow), then $|Tr| \leq 7$ (resp. $|Tr| \leq 13$) and we are done. 

	If $Tr$ is a $\theta$-graph, let $u_1,u_2$ be the vertices of degree 3 in $Tr=\theta_{a,b,c}$ where $u_1,u_2$ are connected by $P_a,P_b,P_c$ respectively. Assume that $a \leq b  \leq c$. If $a+b \geq 10$ or $c \geq 8$, by a similar argument as in Lemma \ref{L56} we can prove that $Tr$ is not $C_5$-saturated, a contradiction. So we have $|T_r| \leq 2+(a-2)+(b-2)+(c-2) \leq 12$.

	By Lemma \ref{L57}, we just need to consider that  $Tr$ is a dumb. We claim that the length of the path connecting two cycles is at most 3. Otherwise let $x_1x_2 \ldots x_p$ be the path with $p \geq 4$ and $x_1,x_p$ are the vertices lie in the two cycles. Note that either $x_1x_3$ or $x_1x_4$ is an admissible nonedge. But $G+x_1x_3$ or $G+x_1x_4$ can only create a $C_3$ or a $C_4$, a contradiction with $Tr$ being $C_5$-saturated. Hence the length of the path connecting two cycles is at most 3. By  Lemma \ref{L56}, we have $|Tr| \leq 15$. \q

\vspace{0.3em}
	
By Corollary \ref{C53}, the component of $Br$ with radius 2 can only be a star  with all vertices of degree one lie in the same part. By Lemmas \ref{L52} and \ref{L55}, we assume that $Br_2,\ldots,Br_m$ are components of radius 1 and $1\le rad(Br_1) \le 2$. Let $V(Br_j)=\{x_j\}$ if $ rad(Br_j)=1$, else $x_1$ be the central vertex  of $Br_1$ if $rad(Br_1) =2$. Then  $N(x_j)\cap V(Tr)=\{\overline{r}_j\}$ for $j\in\{1,\ldots,m\}$. If $rad(Br_1) = 2$, then we will assume $V(Br_1) \subseteq V_1 \cup V_2$ such that all vertices of degree one in $Br_1$ lie in $V_1$ and in this case, $V_2\cap V(Br_1)=\{x_1\}$.
\vspace{0.5em}

\noindent{\bf Claim 2. }(1) $|V(Tr) \cap V_i| \, \leq 8$ for all $i \in [4]$. (2) For all $i \in \{2,3,4\}$, there is $j\in\{2,\ldots,m\}$ such that $V_i\cap V(Br_j)\not=\emptyset$.
\vspace{0.0em}

\noindent{\bf Proof of Claim 2. }By Claim 1, it is easy to check that $|V(Tr) \cap V_i| \, \leq 8$ for all $i \in [4]$ whether $Tr$ is a cycle, a $\theta$-graph, a dumb or a bow.  Since $|V(Br_1) \cap V_2|\, \leq 1$ and $n \geq 10$, the result holds. \q
\vspace{0.3em}

   By Claim 2(2), we assume $x_i \in V_i$ for $i\in \{2,3,4\}$. Then $\overline{r}_i\not=\overline{r}_j$ for distinct $i,j\in \{2,3,4\}$ (otherwise  $G+x_ix_j$ can only form a $C_3$).
By  Claim 2(1) and $n\ge 10$, we have $V_1\setminus V(Tr)\not=\emptyset$. We complete the proof by discuss two cases. 
\vspace{0.3em}	
	
	{\bf Case 1. } There exists $j\in\{j|rad(Br_j) = 1,~1\le j\le m\}$ such that $V_1\cap V(Br_j)\not=\emptyset$.

Assume without loss of generality that $j = 1$, i.e.  $x_1 \in V_1\cap V(Br_1)$. Then $\overline{r}_1\not=\overline{r}_i$ for $i\in \{2,3,4\}$ (otherwise  $G+x_1x_i$ can only form a $C_3$, $i\in \{2,3,4\}$). Since $x_ix_j$ are admissible nonedges in $G$, there exists a $P_3$ connecting $\overline{r}_i$ and $\overline{r}_j$ for all different $i,j \in [4]$. Note that $ |Tr|\leq e(Tr) \leq |Tr|+1$. By Lemma \ref{L57},
 $Tr$  can only be one of the following three types (Figure 20), says $Tr_1,Tr_2,Tr_3$.
	
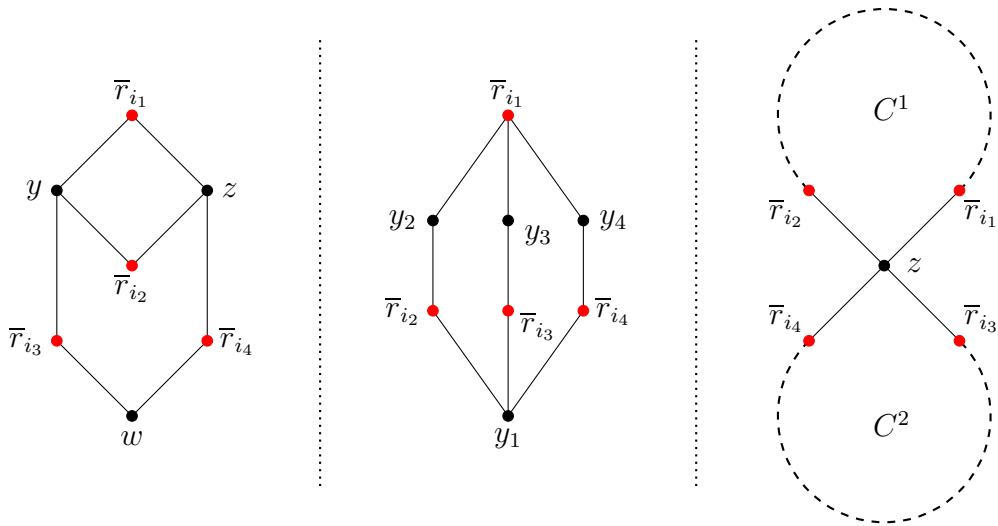
\begin{figure}[H]
\centering

\begin{tikzpicture}[scale=1.0]

\draw  (1,1) -- (1,-1) -- (0,-2) -- (-1,-1) -- (-1,1) -- (0,2) --  (1,1)  (-1,1) -- (0,0) -- (1,1);

\filldraw (-1,1) circle (2pt);
\filldraw[red]  (-1,-1) circle (2pt);
\filldraw (1,1) circle (2pt);
\filldraw[red]  (1,-1) circle (2pt);
\filldraw[red]  (0,2) circle (2pt);
\filldraw[red] (0,0) circle (2pt);
\filldraw (0,-2) circle (2pt);

\draw (1.4,-1) node[align=center]{$\overline{r}_{i_4}$};
\draw (0,2.3) node[align=center]{$\overline{r}_{i_1}$};
\draw (-1.3,1) node[align=center]{$y$};
\draw (1.3,1) node[align=center]{$z$};
\draw (0,-2.3) node[align=center]{$w$};
\draw (-1.4,-1) node[align=center]{$\overline{r}_{i_3}$};
\draw (0,-0.3) node[align=center]{$\overline{r}_{i_2}$};

\draw[dotted,thick] (2.5,3) -- (2.5,-3) ;

\draw (5,2) -- (4,0.6) -- (4,-0.6) -- (5,-2) -- (5,-0.6) -- (5,0.6) -- (5,2) -- (6,0.6) -- (6,-0.6) -- (5,-2) ;

\filldraw[red] (5,2) circle (2pt);
\filldraw[red] (4,-0.6) circle (2pt);
\filldraw[red] (5,-0.6) circle (2pt);
\filldraw[red] (6,-0.6) circle (2pt);
\filldraw (4,0.6) circle (2pt);
\filldraw (5,0.6) circle (2pt);
\filldraw (6,0.6) circle (2pt);
\filldraw (5,-2) circle (2pt);

\draw (5,2.3) node[align=center]{$\overline{r}_{i_1}$};
\draw (3.6,-0.6) node[align=center]{$\overline{r}_{i_2}$};
\draw (5.4,-0.8) node[align=center]{$\overline{r}_{i_3}$};
\draw (6.4,-0.6) node[align=center]{$\overline{r}_{i_4}$};
\draw (3.6,0.6) node[align=center]{$y_2$};
\draw (5.4,0.4) node[align=center]{$y_3$};
\draw (6.4,0.6) node[align=center]{$y_4$};
\draw (5,-2.3) node[align=center]{$y_1$};

\draw[dotted,thick] (7.5,3) -- (7.5,-3) ;

\draw (9,1) -- (10,0) -- (9,-1) (11,1) -- (10,0) -- (11,-1);

\draw[dashed,thick] (11,1) arc(-45:225:1.414213562373095cm and 1.414213562373095cm) ;
\draw[dashed,thick] (9,-1) arc(135:405:1.414213562373095cm and 1.414213562373095cm) ;

\filldraw[red] (9,1) circle (2pt);
\filldraw[red] (11,1) circle (2pt);
\filldraw[red] (9,-1) circle (2pt);
\filldraw[red] (11,-1) circle (2pt);
\filldraw(10,0) circle (2pt);

\draw (11.3,0.7) node[align=center]{$\overline{r}_{i_1}$};
\draw (8.7,0.7) node[align=center]{$\overline{r}_{i_2}$};
\draw (11.3,-0.7) node[align=center]{$\overline{r}_{i_3}$};
\draw (8.7,-0.7) node[align=center]{$\overline{r}_{i_4}$};
\draw (10.4,0) node[align=center]{$z$};
\draw (10.1,2.1) node[align=center]{$C^1$};
\draw (10.1,-2.1) node[align=center]{$C^2$};

\end{tikzpicture}\\

\caption{\centering $Tr_1$ (left), $Tr_2$ (middle), $Tr_3$ (right). The red vertices are $\overline{r}_i$ (under permutation). $Tr_1$, $Tr_2$ are $\theta$-graphs, $Tr_3$ is a bow. All pairs of $\overline{r}_i,\overline{r}_j$ are connected by a $P_3$.}

\end{figure}

	If $Tr=Tr_1$, then $\overline{r}_{i_1}, \overline{r}_{i_2},  w$ (see Figure 20) must lie in the same part by  $Tr$ being $C_5$-saturated. But  $G+x_{i_1}\overline{r}_{i_2}$ would not create a $C_5$, a contradiction. If $Tr=Tr_2$, then $\overline{r}_{i_a}, y_b$ (see Figure 20) must lie in the same part by $Tr$ being $C_5$-saturated for distinct $a,b \in \{2,3,4\}$. Thus $\overline{r}_{i_2}, \overline{r}_{i_3}, \overline{r}_{i_4} , y_2,y_3,y_4$  lie in the same part of $G$, a contradiction.
	
	The case remaining is $Tr=Tr_3$. Let $C^1$ be the cycle containing $z, \overline{r}_{i_1}, \overline{r}_{i_2}$ and $C^2$  the cycle containing $z, \overline{r}_{i_3}, \overline{r}_{i_4}$ (see Figure 20).
 Since $ \overline{r}_{i_2}$ and $ z$ lie in different parts of $G$, either $x_{i_j}\overline{r}_{i_2}$ or $x_{i_j}z$ is an admissible nonedge where $j \in \{1,3,4\}$.
Since $x_i\in V_i$ for $1\le i\le 4$, $\overline{r}_1,\ldots,\overline{r}_4$ are not all lie in one part. We  assume without loss of generality that either $\overline{r}_{i_1} \overline{r}_{i_2}$ or $ \overline{r}_{i_1} \overline{r}_{i_3}$ is an admissible nonedge. If $\overline{r}_{i_1} \overline{r}_{i_2}$ is an admissible nonedge, then $C^1$ is a 6-cycle, say $z \overline{r}_{i_1} y_1 y_2 y_3 \overline{r}_{i_2}z$. Then  $G+x_{i_1}\overline{r}_{i_2}$ or $G+x_{i_1}z$ would only create $C_4, C_6$ or $C_3,C_7$, a contradiction. If $ \overline{r}_{i_1} \overline{r}_{i_3}$ is an admissible nonedge, then either $C^1$ or $C^2$ is a 4-cycle. By symmetry we can assume that $C^1$ is a 4-cycle, say  $z \overline{r}_{i_1} y \overline{r}_{i_2}z$. Then $\overline{r}_{i_1}$ and $\overline{r}_{i_2}$ lie in the same part of $G$ since $Tr$ is $C_5$-saturated. Thus  $G+x_{i_1}\overline{r}_{i_2}$ would only create a $C_4$, a contradiction. 
\vspace{0.3em}	

	{\bf Case 2. }$V_1\setminus V(Tr) \subseteq  V(Br_1)$ and $rad(Br_1)=2$.

Recall we assume all vertices of degree one in $Br_1$ lie in $V_1$ and $x_1\in V_2$. Then $\overline{r}_1\in V_1$.
 Pick $x \in V(Br_1) \cap V_1$. Then $N(x)=\{x_1\}$. Since  $x_i x$ for $2\le i\le 4$ are admissible nonedges  and $G$ is $C_5$-saturated, $\overline{r}_1\overline{r}_i\in E(G)$ for all $i\in \{2,3,4\}$ which implies  $d(\overline{r}_1) \geq 3$ and $\overline{r}_2,\overline{r}_3,\overline{r}_4 \not \in V_1$. Since $x_i \overline{r}_1$  for $2\le i\le 4$ are admissible nonedges, there exists $P_4$ connecting $\overline{r}_1$ and $\overline{r}_i$. Note that $ |Tr|\leq e(Tr) \leq |Tr|+1$. By Lemma \ref{L57}, $Tr$ can only be one of the following two types (Figure 21), says $Tr_4,Tr_5$. Since $x_{i_2}, x_{i_3}$ lie in different parts, either $x_{i_2}\overline{r}_{i_4}$ or $x_{i_3}\overline{r}_{i_4}$ is an admissible nonedge.

\begin{figure}[H]
\centering

\begin{tikzpicture}[scale=1.2]

\draw (1,1) -- (0,2) -- (-1,1) -- (0,0) -- (1,1) (-1,1) -- (0,1) -- (1,1);

\filldraw[blue] (-1,1) circle (2pt);
\filldraw (1,1) circle (2pt);
\filldraw[red]  (0,2) circle (2pt);
\filldraw[red] (0,0) circle (2pt);
\filldraw[red] (0,1) circle (2pt);

\draw (0,2.3) node[align=center]{$\overline{r}_{i_2}$};
\draw (-1.3,1) node[align=center]{$\overline{r}_{1}$};
\draw (1.3,1) node[align=center]{$z$};
\draw (0,-0.3) node[align=center]{$\overline{r}_{i_3}$};
\draw (0,0.7) node[align=center]{$\overline{r}_{i_4}$};

\draw[dotted,thick] (2.5,3) -- (2.5,-1) ;

\draw (6,1) -- (5,2) -- (4,1) -- (5,0) -- (6,1) (8,1) -- (7,2) -- (6,1) -- (7,0) -- (8,1) ;

\filldraw (4,1) circle (2pt);
\filldraw[blue] (6,1) circle (2pt);
\filldraw[red]  (5,2) circle (2pt);
\filldraw[red] (5,0) circle (2pt);

\filldraw (8,1) circle (2pt);
\filldraw[red]  (7,2) circle (2pt);
\filldraw (7,0) circle (2pt);

\draw (5,2.3) node[align=center]{$\overline{r}_{i_2}$};
\draw (6.4,1) node[align=center]{$\overline{r}_{1}$};
\draw (3.7,1) node[align=center]{$z$};
\draw (5,-0.3) node[align=center]{$\overline{r}_{i_3}$};
\draw (7,2.3) node[align=center]{$\overline{r}_{i_4}$};
\draw (8.4,1) node[align=center]{$w$};
\draw (7,-0.3) node[align=center]{$y$};

\end{tikzpicture}\\

\caption{\centering $Tr_4$ (left), $Tr_5$ (right). The red vertices are $\overline{r}_2,\overline{r}_3,\overline{r}_4$ (under permutation), and the blue vertex is $\overline{r}_1$. $Tr_4$ is a $\theta$-graph, $Tr_5$ is a bow. All pairs of $\overline{r}_i\overline{r}_1$ are connected by a $P_4$.}

\end{figure}
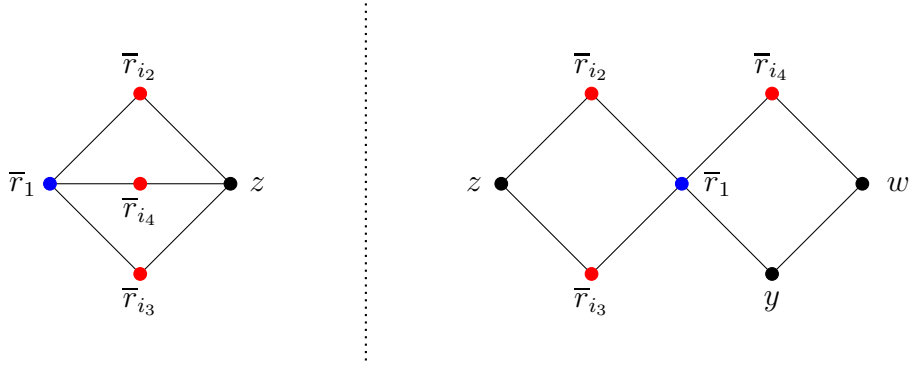

	If $Tr=Tr_4$, then  $G+x_{i_2}\overline{r}_{i_4}$ or $G+x_{i_3}\overline{r}_{i_4}$ would only create $C_4$ or $C_6$, a contradiction. If $Tr=Tr_5$, then  $G+x_{i_2}\overline{r}_{i_4}$ or $G+x_{i_3}\overline{r}_{i_4}$ would only create $C_4$ or $C_6$, a contradiction.

Thus  we finished our proof of Theorem \ref{T17}. \qed	
\vspace{1em}

\noindent{\bf Proof of Theorem \ref{T16}. }Let $G \in Sat(  K_3^n, C_5)$ and $G=V_1\cup V_2\cup V_3$. Assume $Tr$ is the trunk of $G$ and $Br$ is the branch of $G$ with $m$  components $Br_1,\ldots,Br_m$ and $N(V(Br_j))=\{\overline{r}_j\}$ for $1\le i\le m$. By Theorem \ref{T12}, $e(G) \geq 3n=|G|$.

If $n=2$,  $C_6=x_1x_2x_3x_4x_5x_6$ with $x_i, x_{i+3} \in V_i$ for $i \in [3]$ is a $C_5$-saturated tripartite graph. Hance $e(G)\le 6$ and we are done.

	For $n \geq 3$, by the construction in Figure 5, $e(W_*^{(5,3,n)})= 3n+1$. We first check that $W_*^{(5,3,n)}$ is $C_5$-saturated. If $W_*^{(5,3,n)}$ has a $C_5$, then it must lie in $Tr= G[A]$ where $A=\{u^*,w,u_3,u_2\}$. Since $|Tr|=4$, we have $W_*^{(5,3,n)}$ is $C_5$-free.

	Let $v_i \in B_i$ and $N(v_i)=\{u_i\}$ (see Figure 5). The pairs of vertices which are the endpoints of admissible nonedge in $W_*^{(5,3,n)}$ are all connected by a $P_5$\,:
\begin{table}[H]
\centering
\begin{tabular}{|c|c|}
\hline
\quad Admissible nonedge \quad & \quad The $P_5$ connecting them \quad  \\
\hline
\, $v_iv_j$ \, ($1 \leq i < j \leq 3$) &  $v_iu_iu^*u_jv_j$ \\
\hline
\, $v_1 u_i$ \, ($2 \leq i \leq 3$) &  $v_1u_1u^*u_{5-i}u_i$ \\
\hline
\, $v_i z$ \, ($2 \leq i \leq 3$, $ z \in \{u^*,w\}$) & \, $v_iu_iz^*u_{5-i}z$ \, ($ z^* \in \{u^*,w\} \backslash \{z \}$) \\
\hline
\, $v_3 u_1$ \,  & \, $v_3u_3u_2u^*u_1$ \,  \\
\hline
\, $u_1w$,\, $u_1u_2$  & \, $u_1u^*u_2u_3w$,\,  $u_1u^*u_3wu_2$  \\
\hline

\end{tabular}
\end{table}
\noindent Thus $W_*^{(5,3,n)}$ is $C_5$-saturated and we have $sat(K_3^n,C_5) \leq 3n+1$.

	 Suppose that $e(G)=3n=|G|$. By Observation \ref{O51} (iii), $|Br|=|G|-|Tr|$ and $|G|=3n$, we have $e(Tr)= |Tr|$. By Observation \ref{O51} (i) and (ii),
 $Tr$ is a cycle. By Lemma \ref{L56}, $|Tr| \leq 7$.
\vspace{0.5em}	

\noindent{\bf Claim. }$|Tr| = 6$.
\vspace{0.0em}

\noindent{\bf Proof of Claim. }Suppose $Tr=u_1u_2u_3u_1$ with $u_i \in V_i$. Assume $x \in V(Br)\cap N(u_2)\cap V_1$. Then  $G+xu_3$ would not have a $C_5$, a contradiction. Suppose $Tr= u_1u_2u_3u_4u_1$. Since $Tr$ is $C_5$-saturated, we have $u_i,u_{i+2}$ lie in the same part. Assume that $u_i \in V_i $ for $i \in [2]$. By symmetry we can assume that $u_1$ has a neighbour $x \in V(Br)$. Then  $G+xu_3$ would not have a $C_5$, a contradiction. Since $G$ is $C_5$-saturated, $|Tr| \neq 5$. Suppose $Tr= u_1 \ldots u_7u_1$. Since $Tr$ is $C_5$-saturated, we have $u_i,u_{i+2}$ lie in the same part. But then $u_1,u_7$ lie in one part of $G$, a contradiction. Thus $Tr$ must be a $C_6$. \q
\vspace{0.3em}

Let $Tr= u_1 \ldots u_6u_1$. Since $Tr$ is $C_5$-saturated, we have $u_i,u_{i+3}$ lie in the same part of $G$. Assume that $u_i \in V_i$ for $i \in [3]$. Since $n \geq 3$, by Lemma \ref{L55}, there exists a $Br_j$ such that $rad(Br_j)=1$. Assume $Br_j=\{x\}$, $x \in V_1$ and $xu_2\in E(G)$. Then  $G+xu_3$ would not create a $C_5$, a contradiction. Hence $e(G) \geq 3n+1$, and we finished the proof of Theorem \ref{T16}. \qed

\section{Concluding remark}

Determining the value of $sat(K_k^n, C_{\ell})$ is an  interesting problem in partite saturation number. Although we have determined the value for $\ell > k \geq 3$ and $\ell =2r \geq 6, k=2$ asymptotically except for $(\ell,k)=(4,4)$, there is still the other half of the problem that needs to be solved urgently. We give our conjectures on such problems.

\begin{prob}\label{P61}
Determine the value of $sat(K_k^n, C_{\ell})$ for all $ \ell \geq 4, k \geq 3$ and even $\ell \geq 8, k=2$.
\end{prob}

\begin{conj}\label{C62}
For $\ell \geq 4$ and $n_1,n_2 \geq \ell+2$, $sat(K_{n_1,n_2},C_{2\ell})=n_1+n_2+\ell^2-3\ell+1$.
\end{conj}

\begin{conj}\label{C63}
For $n \geq 1$, $sat(K_4^n,C_4)=5n-1$.
\end{conj}

%
	
	For a given graph family $\mathcal{F}$, a graph $H$ is said to be $\mathcal{F}$-saturated in $G$ if $H$ does not contain any copy of all graphs in $ \mathcal{F}$, but the addition of a nonedge in $E(G) \setminus E(H)$ would form a copy of some $F \in \mathcal{F}$. Similarly we can define saturation number $sat(G, \mathcal{F})$ and partite saturation number $sat(K_k^n,\mathcal{F})$. Let $\mathcal{C}_{\geq \ell }$ be the graph families of cycles with length at least $\ell$, i.e. $\mathcal{C}_{\geq \ell } = \{ C_{\ell},C_{\ell+1}, \ldots \}$. Although Theorem \ref{T110} shows that the Construction I is not optimal, we believe that it is tight when we consider $\mathcal{C}_{\geq \ell}$-saturation. So we  have following conjecture.
	
\begin{conj}\label{C64}
For $\ell >k \geq 3$, $\ell \geq 6$ and sufficiently large $n$,
\begin{equation*}
sat(K_k^n,C_{\geq \ell})\; =
\; kn- \ell+1 + \lfloor \frac{\ell-2}{2} \rfloor ^2  + 2 \left( \ell-1-2 \, \lfloor \frac{\ell-2}{2} \rfloor \right) \lfloor \frac{\ell-2}{2} 
\rfloor \, ,
\end{equation*}
and the extremal graphs are all in $\Omega^{(\ell,k,n)}$.
\end{conj}

\section*{Acknowledgement}
This paper is supported by the National Natural Science Foundation of China (No.~12401445, 12171272, 12161141003); and by Beijing Natural Science Foundation (No.~1244047), China Postdoctoral Science Foundation (No.~2023M740207).

\end{document}